\documentclass[12pt]{amsart}
\usepackage{hyperref}

\hypersetup{
  colorlinks= true, 
  urlcolor   = red, 
  linkcolor  = blue, 
  citecolor  = blue 
}

\newcommand{\modd}[1]{\,(\text{mod}\, {#1})}
\def\Re{\text{Re}\,}
\def\Im{\text{Im}\,}

\newtheorem{theorem}{Theorem}
\newtheorem{lemma}{Lemma}
\theoremstyle{definition}
\newtheorem{example}{Example}

\makeatletter
\def\specialsection{\@startsection{section}{1}%
  \z@{\linespacing\@plus\linespacing}{.5\linespacing}%
  {\normalfont}}
\def\section{\@startsection{section}{1}%
  \z@{.7\linespacing\@plus\linespacing}{.5\linespacing}%
  {\normalfont\scshape}}
\makeatother

\title[Multiplicative Functions and Maass Forms]{Harmonic Analysis on the Positive Rationals II: Multiplicative Functions and Maass Forms}

\author{Peter D. T. A. Elliott}
\address{Department of Mathematics, University of Colorado Boulder, Boulder, Colorado 80309-0395 USA}
\email{pdtae@euclid.colorado.edu}

\author{Jonathan Kish}
\address{Department of Mathematics, University of Colorado Boulder, Boulder, Colorado 80309-0526 USA}
\email{jonathan.kish@colorado.edu}

\begin{document}
\parskip10pt
\parindent15pt
\baselineskip15pt

\maketitle
\setcounter{lemma}{18}


\section{Introduction}\label{JP_elliott_04_sec_automorphic_forms}


Ramanujan's arithmetic $\tau$-function is defined by
\begin{equation}
\sum\limits_{n=1}^\infty \tau(n) x^n = x \prod\limits_{j=1}^\infty (1-x^j)^{24}, \quad |x| < 1. \notag
\end{equation}
Thirty years ago the first author conjectured that it satisfied the Central Limit condition
\begin{equation}
\mathop{\phantom{1}x^{-1}\sum 1 \phantom{x^{-1}}  }\limits_{\substack{n \le x \\ |\tau(n)| \le n^{11/2} \left((\log x)^{-1/2} \exp\left(z c(\log \log x)^{1/2} \right)\right)}} \to \frac{1}{\sqrt{2\pi}} \int_{-\infty}^z e^{-u^2/2} \, du, \quad x \to \infty. \notag
\end{equation}
This is the weak convergence of measures on the real line and, in view of the continuity of the limiting distribution, would be uniform in real $z$.  An appropriate value for the positive constant $c$ would be $\left(\pi^2/12 + 1/2\right)^{1/2}$, c.f. Elliott \cite{elliott1981MFramanujan}, \cite{elliott2012centrallimitramanujantau}, \cite{elliott2012CLTeigenforms}, \cite{elliott2013CLTcuspformsbasilgordon}.

Implicit in this conjecture is the conjecture of Lehmer \cite{lehmer1943ramanujanstau}, that $\tau$ never vanishes, for otherwise the multiplicativity of $\tau$ would force the limiting measure to have an atom at $z=0$.

The first author proved recently that if
\begin{equation}
\sum\limits_{n=1}^\infty a_n e^{2\pi i nz}, \quad a_1 = 1, \ \Im(z)>0, \notag
\end{equation}
is an elliptic holomorphic new form of weight $k \ge 2$, level $N$ and nebentypus $\psi: (\mathbb{Z}/N\mathbb{Z})^* \to \mathbb{C}^*$, then for each real $z$
\newlength{\Myl}
\settowidth{\Myl}{$|a_n|^2 n^{-(k-1)}$}
\begin{equation}
\mathop{ \makebox[\Myl][r]{$(\alpha x)^{-1}$} \sum |a_n|^2 n^{-(k-1)} } \limits_{\substack{n \le x \\ |a_n| \le n^{(k-1)/2} \exp\left(A(x) + zB(x)\right)}} \to \frac{1}{\sqrt{2\pi}} \int_{-\infty}^z e^{-u^2/2} \, du, \quad x \to \infty,  \notag
\end{equation}
where
\begin{equation}
A(x)  = \sum\limits_{p \le x} |a_p|^2 p^{-k} \log\left(|a_p| p^{-(k-1)/2} \right), \notag
\end{equation}
\begin{equation}
B(x)  = \left( \sum\limits_{p \le x} |a_p|^2 p^{-k} \left( \log \left(|a_p| p^{-(k-1)/2} \right) \right)^2 \right)^{1/2} \ge 0, \notag
\end{equation}
the sums taken over primes for which $a_p \ne 0$, and
\begin{equation}
\alpha = \lim\limits_{x \to \infty} x^{-1} \sum\limits_{n \le x} |a_n|^2 n^{-(k-1)}, \notag
\end{equation}
c.f. \cite{elliott2013CLTcuspformsbasilgordon}.

For those forms that do not possess complex multiplication in the sense of Ribet \cite{ribet1977galoiseigenforms}, i.e. whose coefficients $a_p$ do not vanish identically on the inertial primes of an imaginary quadratic extension of the rationals, which includes the form with $k=12$, $N=1$, corresponding to Ramanujan's function, the version of the Sato-Tate conjecture established by Barnet-Lamb, Geraghty, Harris and Taylor \cite{barnetlambgeraghtyharristaylor2011satotate}, ensures that we may replace $B(x)$ by $\left((\pi^2/12 + 5/8)\log\log x\right)^{1/2}$ and estimate $A(x)$ to be $\left(1/4 + o(1)\right)\log\log x$ as $x \to \infty$.

The basic result remains valid whether or not the cusp form under consideration has complex multiplication.

If we identify the above cusp forms with functions in the Hilbert space\\
$L^2( \Gamma_0(N) \backslash SL(2,\mathbb{R}) )$ with respect to the appropriate Haar measure, c.f. Gelbart \cite{gelbart1975automorphicadele}, then the introduction of the weights $|a_n|^2 n^{-(k-1)}$ appears natural and subdues the effect of vanishing coefficients.  In order to remove the weights in favour of $1$, the vanishing of the $a_n$ has to be met, since it directly affects the formulation of appropriate limit laws.

Appealing to an $\ell$-adic representation of Deligne, and the Chebotarev density theorem, Serre \cite{serre1981chebotarev}, 1981, showed that for each member of this same class of holomorphic cusp forms, once again without complex multiplication, there is a positive constant $b$ so that amongst the coefficients $a_p$, $2 \le p \le x$, at most $O\left(x(\log x)^{-1-b}\right)$ vanish.  As a consequence the density
\begin{equation}
\rho = \lim\limits_{x \to \infty} x^{-1} \sum\limits_{\substack{n \le x \\ a_n \ne 0}} 1 \notag
\end{equation}
exists and is positive.

This may be applied to the forms of weight $2$ that correspond to elliptic curves, the existence of which was conjectured by Taniyama and Shimura, and established by Wiles \cite{wiles1995fermatannals}, Taylor and Wiles \cite{taylorwiles1995hecke}, Breuil, Conrad, Diamond, Taylor \cite{breuilconraddiamondtaylor2001wild3adic}.  By a result of Elkies \cite{elkies1987supersingular}, in these cases  $\rho$ is always less than $1$.

Otherwise for each form with complex multiplication there is a positive constant $A$ and an asymptotic estimate
\begin{equation}
\sum\limits_{\substack{n \le x \\ a_n \ne 0}} 1 \sim Ax(\log x)^{-1/2}, \quad x \to \infty, \notag
\end{equation}
Serre, loc. cit.

To obviate appeal to an hypothesis of Lehmer type, a natural analogue, for elliptic holomorphic new forms without complex multiplication, of the conjecture concerning Ramanujan's function would be that
\begin{equation}
\mathop{ \phantom{1}(\rho x)^{-1} \sum 1 \phantom{(\rho x)^{-1}} } \limits_{\substack{n \le x \\ 0 < |a_n| n^{-(k-1)/2} \le (\log x)^{-1/2} \exp\left(z\delta(\log\log x)^{1/2} \right)}} \to \frac{1}{\sqrt{2\pi}} \int_{-\infty}^z e^{-u^2/2} \, du, \quad x \to \infty, \notag
\end{equation}
where $\rho$ is the density of nonvanishing of the coefficients, $\delta^2 = \pi^2/12 + 1/2$.

An improvement of the Sato-Tate asymptotic estimate of Barnet-Lamb, Geraghty, Harris and Taylor loc. cit., to within an error of $o\left(x(\log x)^{-3}\right)$ for primes in the interval $2 \le p \le x$, would suffice to establish the general conjecture.

For forms with complex multiplication the corresponding frequency to consider would be
\begin{equation}
 \mathop{ \phantom{1}\left(Ax(\log x)^{-1/2} \right)^{-1} \sum 1 \phantom{\left(Ax(\log x)^{-1/2} \right)^{-1}}  } \limits_{\substack{n \le x \\ 0 < |a_n| n^{-(k-1)/2} \le (\log x)^{-1/2} \exp\left(z\lambda(\log\log x)^{1/2} \right)}} \notag
\end{equation}
with $\lambda^2 = \pi^2/24$, the desired result being in principle already within reach.

A somewhat more elaborate discussion of these matters, together with examples, may be found in Elliott \cite{elliott2013CLTcuspformsbasilgordon}.

Given a positive integer $D$, introduction of Dirichlet characters would allow the method of Serre to be extended to the consideration of coefficients of a form on an arithmetic progression $\modd{D}$.  Less clear is how uniform the various estimates might be in $D$.

From the viewpoint of functional analysis or, more widely, of group representations, it is natural to expect analogues of the above results for the coefficients of other automorphic cusp forms that are eigenfunctions of the appropriate Hecke operators, in particular for the coefficients of Maass wave forms.

Let
\begin{equation}
f = \sum\limits_{n \ne 0} a_n (2\pi y)^{1/2} K_s(2\pi |n| y) e^{2\pi i n x} \notag
\end{equation}
be a nonzero Maass cusp form, attached to the action of $SL(2,\mathbb{Z})$ on the complex upper half-plane, that is an eigenfunction of the appropriate Hecke operators and renormalised so that $a_n$ is a multiplicative function of $n$.  An analogue of the Ramanujan-Petersson conjecture $|a_p| \le 2p^{(k-1)/2}$ for holomorphic forms, established by Deligne, is unavailable, even less an analogue of the Sato-Tate conjecture.

In the present paper we address the question of how often the coefficients of a Maass form vanish on a given arithmetic progression or, when real, have a particular sign.  In order that the argument have the widest possible application we assume as little as possible of the analytic properties of the $L$-functions attached to the forms, this being anyway currently forced upon us.

The above example of a Maass form may be realised in terms of an irreducible component attached to a representation of the group $GL_2(A_\mathbb{Q})$, where $A_F$ denotes the ring of ad\`eles of the algebraic extension  $F$ of the rationals.

In a series of works from Gelbart and Jacquet \cite{gelbartjacquet1978automorphicrepresentations}, Kim and Shahidi \cite{kimshahidi2000functorialcube}, Kim \cite{kim2003functorialitysymmetricsquare}, amongst others, serious strides have been made in Langland's functoriality program.  With specific exceptions characterised in Kim and Shahidi \cite{kimshahidi2002cuspidality}, the second, third and fourth symmetric product $L$-functions attached to a normalised nonzero cusp form, derived from a representation of $GL_2(A_F)$ and an eigenfunction of the appropriate Hecke operators, may be realised as $L$-functions of a cusp form derived from representations of $GL_k(A_F)$, $k = 3$, $4$, $5$, respectively.  In particular, they have analytic continuations over the whole complex plane.  Moreover, for immediate application to Analytic Number Theory, adequate control is available over symmetric product $L$-functions of order up to $9$, c.f. Kim and Shahidi loc. cit.

The Fourier coefficients $a_p$ at the infinite cusp of the Maass form $f$ over $\mathbb{Q}$ are real and there is an Euler product representation for the attached $L$-function:

\begin{align}
\sum\limits_{n=1}^\infty a_n n^{-s} & = \prod\limits_p \left( 1 - a_p p^{-s} + p^{-2s} \right)^{-1} \notag \\
& = \prod\limits_p \left(1-\alpha_p p^{-s} \right)^{-1} \left(1-\beta_p p^{-s} \right)^{-1}, \quad \Re(s) > 1, \notag
\end{align}
where, courtesy of Kim and Shahidi loc. cit., $\max\left(|\alpha_p|, |\beta_p| \right) \le p^{1/9}$.  A direct calculation with Euler products shows that
\begin{equation}
\exp\left( \sum\limits_p a_p^4 p^{-s} + h(s) \right) = L(\text{sym}^2 \times \text{sym}^2, s) L(\text{sym}^2, s)^2\zeta(s), \notag
\end{equation}
where
\begin{equation}
h(s) = \sum\limits_p \sum\limits_{k=2}^\infty k^{-1}\left( \alpha_p^k + \beta_p^k \right)^4 p^{-ks}, \notag
\end{equation}
the first $L$-function is the Rankin-Selberg product of the symmetric square $L$-function attached to (the Maass form) $f$ realised, under the Gelbart-Jacquet lift, as the $L$-function of a Maass form representing the action of $SL(3,\mathbb{Z})$ on the generalised upper half-plane $GL(3,\mathbb{R})/ \left(O(3,\mathbb{R}) \cdot \mathbb{Z}_3\right)$ and so on; c.f. Goldfeld \cite{Goldfeld2006}.

Since $\alpha^k + \beta^k$ is a symmetric function of $\alpha$, $\beta$ hence a polynomial in $\alpha + \beta$, $\alpha \beta$, it is real.  The function $h(s)$ is analytic in $\Re(s) > 17/18$ and in the half-plane $\Re(s) \ge 26/27$ satisfies the bound
\begin{equation}
|h(s)| \le 8 \sum\limits_p \sum\limits_{k=2}^\infty p^{-k(\sigma - 4/9)} \le 32 \zeta(28/27). \notag
\end{equation}
The Rankin-Selberg product and the Riemann zeta function each have a simple pole at $s=1$, otherwise the above $L$-functions are analytic in an open set that contains the half-plane $\Re(s) \ge 1$.

Following classical methodology, an application of the Wiener-Ikehara theorem provides the asymptotic estimate
\begin{equation}
\sum\limits_{p \le x} a_p^4 \log p \sim 2x, \quad x \to \infty. \notag
\end{equation}

A like consideration of the Rankin-Selberg product of the standard $L$-function attached to $f$ yields the corresponding
\begin{equation}
\sum\limits_{p \le x} a_p^2 \log p \sim x, \quad x \to \infty. \notag
\end{equation}

For this particular example more is known, but we shall often not need it.  Suffice it, for the moment, that after an application of the Cauchy-Schwarz inequality we obtain the lower bound
\begin{equation}
\sum\limits_{\substack{p \le x \\ a_p \ne 0}} \log p \ge \left(1/2 + o(1)\right)x, \quad x \to \infty. \notag
\end{equation}

Indeed, we shall often only assume that for positive constants $c_1$ and $c_2$,
\begin{equation}
\sum\limits_{\substack{p \le x \\ a_p \ne 0}} \log p \ge c_1x, \quad x \ge c_2; \notag
\end{equation}
and that
%
%
\begin{equation}
\sum\limits_{\substack{u < p \le v \\ a_p \ne 0}} p^{-1} \ge (1/2) \sum\limits_{u < p \le v} p^{-1} + O(1), \notag
\end{equation}
uniformly in $2 \le u \le v$, each of which, after Heilbronn and Landau \cite{heilbronnlandau1933bochner, heilbronnlandau1933anwendungenwienerschen}, may be deduced from the analyticity and nonvanishing of $\zeta(s)^{-1} \sum_{n=1}^\infty a_n^2 n^{-s}$ and $\zeta(s)^{-2} \sum_{n=1}^\infty a_n^4 n^{-s}$ in a proper disc $|s-1| < c_3$.

In general, for the purposes of establishing central limit theorems involving mean-square coefficients of Maass forms it currently suffices to possess first analytic properties of appropriate symmetric $L$-functions of order up to the sixth.

As will be shown in the present paper, when counting the frequency with which Maass coefficients on a given residue class do not vanish, or have a given sign, it suffices to have such properties for orders up to the fourth.

We address these questions in the final section.  With so little to hand almost all the weight of argument will be borne by the multiplicativity of an appropriately defined arithmetic function.  Since classical methodology with its appeal to analytic continuation is unavailable, the argument will operate under a completely different aesthetic.

To ensure a reasonably self-contained presentation, the account of Theorem \ref{JP_elliott_04_thm_02} given here occasionally overlaps that of Theorem 1 in the previous paper, I.  The emphasis is, however, quite different.


\section{Multiplicative Functions}\label{JP_elliott_04_sec_multiplicative_functions}


A systematic study of multiplicative functions with values in the complex unit disc, initiated by Delange in 1961 \cite{delange1961surlesfonctions}, received strong impulses from Wirsing \cite{wirsing1967}, 1967, and Hal\'asz \cite{halasz1968}, 1968.  Implicitly the multiplicative function was largely compared to the function that is identically 1 on the positive integers.

In the present section we begin in general terms and gradually narrow the focus towards the problem raised in the previous section.  We consider the mean-value of a multiplicative function that may assume large values but whose support on the primes may be seriously decreased, rendering it far from the function 1.

To gain perspective on matters of size the next three results concern multiplicative functions with values in a real interval, $[0,\beta]$.


\makeatletter{}

\begin{lemma} \label{elliott_duality_lem_02_02}
The estimate
\begin{equation}
\sum\limits_{2 \le n \le x} g(n) \le \left(\frac{x}{\log x} + \frac{10x}{(\log x)^2} \right) \Delta \sum\limits_{n \le x} \frac{g(n)}{n} \notag
\end{equation}
with
\begin{equation}
\Delta = \sup\limits_{1 \le y \le x} y^{-1} \sum\limits_{q \le y} g(q)\log q, \notag
\end{equation}
the sum over prime-powers $q$, holds uniformly for all real nonnegative multiplicative functions $g$, and all $x \ge 2$.  Moreover,
\begin{equation}
\frac{1}{\log x} \sum\limits_{n \le x} \frac{g(n)}{n} \ll \exp\left(\sum\limits_{q \le x} \frac{g(q)-1}{q} \right) \notag
\end{equation}
with an absolute implied constant.  In particular, if $g(q) \le M$ on the prime-powers, then
\begin{equation}
x^{-1} \sum\limits_{2 \le n \le x} g(n) \ll M \exp\left(\sum\limits_{q \le x} \frac{g(q)-1}{q} \right). \notag
\end{equation}
\end{lemma}    

\noindent \emph{Proof}.  A complete proof may be found in \cite{elliott1997}, Lemma 2.2.

\makeatletter{}

\begin{lemma}\label{JP_elliott_04_lem_02}
Under the conditions of Lemma \ref{elliott_duality_lem_02_02} assume that $g(p) \le \beta$ uniformly on the primes $p$ and that the series $\sum g(q) q^{-1}$, taken over the prime-powers $q=p^k$ with $k \ge 2$, converges.  Then
\begin{equation}
1 \ll \prod\limits_{p \le x} \left(1+g(p)p^{-1} \right)^{-1} \sum\limits_{n \le x} g(n)n^{-1} \ll 1, \quad x \ge 1. \notag
\end{equation}
\end{lemma}     

\makeatletter{}

\begin{lemma}\label{JP_elliott_04_lem_03}
Assume, further, that
\begin{equation}
\sum\limits_{p \le x} g(p) \log p \ge cx \notag
\end{equation}
for a positive $c$ and all sufficiently large $x$.  Then
\begin{equation}
\sum\limits_{n \le x} g(n) \gg x \exp\left(\sum\limits_{p \le x} -\left(1-g(p)\right)p^{-1} \right), \quad x \ge 1. \notag
\end{equation}
\end{lemma}     

\makeatletter{}

\noindent \emph{Proofs}.  Lemmas \ref{JP_elliott_04_lem_02} and \ref{JP_elliott_04_lem_03} will be established in two passes, the first with $\beta = 1$.

The upper bound in Lemma \ref{JP_elliott_04_lem_02} follows from the inequality
\begin{equation}
\sum\limits_{n \le x} g(n)n^{-1} \le \prod\limits_{p \le x} \left(1+\sum\limits_{k=1}^\infty g(p^k)p^{-k}\right) \notag
\end{equation}
and the convergence of the series $\sum_{p,k \ge 2} g(p^k)p^{-k}$.

Towards the lower bound we may assume that $g(p)=0$ if $x^\varepsilon < p \le x$, $0 < \varepsilon \le 1/2$, and define the multiplicative function $h$ by $h(p) = 1-g(p)$, $h(p^k)=0$, $k \ge 2$.  On squarefree integers the Dirichlet convolution $h * g$ assumes the value 1.  Restricting summation to squarefree integers,
\begin{align}
\frac{6}{\pi^2}x & + O(x^{1/2}) = \sum\limits_{ab \le x} h(a)g(b) \notag \\
& \le \sum\limits_{a \le x^\varepsilon} h(a) \sum\limits_{b \le x/a} g(b) + \sum\limits_{b \le x^{1-\varepsilon}} g(b) \sum\limits_{a \le x/b} h(a) \notag \\
& = \Sigma_1 + \Sigma_2. \notag
\end{align}

By Lemma \ref{elliott_duality_lem_02_02},
\begin{equation}
\sum\limits_{b \le x/a} g(b) \ll xa^{-1} \exp\left(-\sum\limits_{p \le x/a} h(p)p^{-1} \right), \notag
\end{equation}
from which
\begin{equation}
\Sigma_1 \ll x\sum\limits_{a \le x^\varepsilon} h(a)a^{-1} \exp\left(-\sum\limits_{p \le x^\varepsilon} h(p)p^{-1} - \sum\limits_{x^\varepsilon < p  \le x^{1-\varepsilon}} p^{-1} \right) \ll \varepsilon x, \notag
\end{equation}
since $(x^\varepsilon, x^{1-\varepsilon}]$ is contained in $(x^\varepsilon, x/a]$.

With $\varepsilon$ fixed at a suitably small value, independent of $g$,
\begin{align}
x \ll \Sigma_2 & \ll x \sum\limits_{b \le x^{1-\varepsilon}} g(b) b^{-1} \exp\left( -\sum\limits_{p \le x/b} g(p)p^{-1} \right) \notag \\
& \ll \exp\left( -\sum\limits_{p \le x^\varepsilon} g(p)p^{-1} \right) x \sum\limits_{b \le x} g(b)b^{-1}. \notag
\end{align}
The sum $\sum_{x^\varepsilon < p \le x} p^{-1}$ is bounded in terms of $\varepsilon$ alone and the desired result follows readily.

Under the further assumption of Lemma \ref{JP_elliott_04_lem_03}, once again restricting to squarefree integers, for $x$ sufficiently large
\begin{align}
\sum\limits_{n \le x} & g(n)\log n = \sum\limits_{n \le x} g(n) \sum\limits_{p \mid n} \log p \notag \\
& \ge \sum\limits_{m \le x^{1/2}} g(m) \left(\sum\limits_{p \le x/m} g(p)\log p - \log m \right) \notag \\
& \gg x \sum\limits_{m \le x^{1/2}} g(m)m^{-1} \gg x \prod\limits_{p \le x^{1/2}} \left(1+g(p)p^{-1}\right). \notag
\end{align}
Moreover, the factor $\log n$ does not exceed $\log x$ which, in turn, lies between constant multiples of $\exp\left(\sum_{p \le x} p^{-1} \right)$.

For the second pass we assume that $\beta > 1$.  Lemma \ref{elliott_duality_lem_02_02} again supplies an upper bound, which it is convenient to express in the form
\begin{equation}
\sum\limits_{n \le x} g(n) \ll x(\log x)^{\beta-1} \exp\left(-\sum\limits_{p \le x} (\beta - g(p))p^{-1}\right), \quad x \ge 2. \notag
\end{equation}

We may then follow the argument for the case $\beta = 1$, modifying the function $h$ to satisfy $h(p) = \beta - g(p)$ so that on squarefree integers $(h*g)(n) = \beta^{\omega(n)}$, where $\omega(n)$ denotes the number of distinct primes divisors of $n$.

The bound
\begin{equation}
\sum\limits_{n \le x} \mu(n)\beta^{\omega(n)} \gg x(\log x)^{\beta-1}, \notag
\end{equation}
sufficient to complete the proof of the general version of Lemma \ref{JP_elliott_04_lem_02}, may be derived from the case $\beta=1$ of Lemma \ref{JP_elliott_04_lem_03} via the Dirichlet convolution representation
\begin{equation}
\mu^2\beta^\omega = \mu^2 * \mu^2 * \cdots * \mu^2 * (\beta - r)^\omega \mu^2 * t, \notag
\end{equation}
where the $r$ convolutions of the characteristic function of the squarefree integers are chosen so that $0 < \beta-r \le 1$ and the function $t$, given by
\begin{equation}
\sum\limits_{n=1}^\infty t(n)n^{-s} = \exp\left(\sum\limits_p \sum\limits_{k \ge 2} \left(r+(\beta-r)(-1)^{k+1}\right)\left(kp^{ks}\right)^{-1}\right), \notag
\end{equation}
assumes nonnegative values.

In particular,
\begin{equation}
\sum\limits_{n \le x} \mu(n)^2 \beta^{\omega(n)} \ge \sum\limits_{\substack{v_j \le x^{1/(2r)} \\ j=1, \dots, r}} \mu(v_1)^2 \cdots \mu(v_r)^2 \sum\limits_{n \le x} \mu(n)^2 (\beta-r)^{\omega(n)} \notag
\end{equation}
\begin{equation}
\gg x(\log x)^{\beta - r - 1} \left(\sum\limits_{v \le x^{1/(2r)}} \mu(v)^2 v^{-1} \right)^r \gg x(\log x)^{\beta-1}, \notag
\end{equation}
the last step by the first-pass version of Lemma \ref{JP_elliott_04_lem_02}.

This completes the proof of the general Lemma \ref{JP_elliott_04_lem_02}; which result may be employed to complete the proof of the general Lemma \ref{JP_elliott_04_lem_03}.


Introducing Dirichlet characters $\chi \modd{D}$, we consider a typical innersum in the representation
\begin{equation}
\sum\limits_{\substack{n \le x \\ n \equiv a \modd{D}}} g(n) = \varphi(D)^{-1} \sum\limits_{\chi \modd{D}} \overline{\chi}(a) \sum\limits_{n \le x} g(n)\chi(n), \quad (a,D)=1. \notag
\end{equation}
Under the conditions of Lemma \ref{JP_elliott_04_lem_03}, the innersum corresponding to the principal character $\modd{D}$ lies between multiples of
\begin{equation}
x(\log x)^{-1} \prod\limits_{\substack{p \le x \\ (p,D)=1}} \left(1+g(p)p^{-1}\right), \notag
\end{equation}
comparable to the sum $\sum_{n \le x} g(n)$ itself.

Provided $D$ does not exceed a certain power of $x$, depending upon the value of $c$ in the hypotheses of Lemma \ref{JP_elliott_04_lem_03}, this lower bound remains valid uniformly in $D$.

If $g$ is a multiplicative function with values in the complex unit disc and otherwise unconstrained, then the method of Hal\'asz delivers an estimate
\begin{equation}
x^{-1} \sum\limits_{n \le x} g(n) \ll me^{-m} + T^{-1/2}, \quad x \ge 2, \ T \ge 2, \notag
\end{equation}
where
\begin{equation}
m = \min\limits_{|t| \le T} \sum\limits_{p \le x} \left(1-\Re g(p)p^{it} \right)p^{-1}, \notag
\end{equation}
c.f. Tenenbaum \cite{tenenbaum1995}, III.4 Corollary 6.3, in particular exercise 6 of the same chapter; see, also, \cite{Elliott1979}, Chapter 6.  The leading term in the upper bound may be effectively rearranged as
\begin{equation}
(\log x)^{-1} \prod\limits_{p \le x} \left(1+|g(p)|p^{-1} \right)me^{-\rho} \notag
\end{equation}
with
\begin{equation}
\rho = \min\limits_{|t| \le T} \sum\limits_{p \le x} \left(|g(p)| - \Re g(p)p^{it}\right)p^{-1}. \notag
\end{equation}
Supposing that $\sum_{p \le x} |g(p)|\log p \sim cx$, $0 < c < 1$, $x \to \infty$, it is feasible that $\rho$ becomes unbounded with $x$ but remains $O\left((\log\log\log x )^{1/2}\right)$ in size; an example is constructed following the proof of Theorem \ref{JP_elliott_04_thm_02}.  The saving by the factor $e^{-\rho}$ is then vitiated by the factor $m$ of size $\left(1-c+o(1)\right)\log\log x$.

The following result addresses this phenomenon.


\setcounter{theorem}{1}
\makeatletter{}

\begin{theorem}\label{JP_elliott_04_thm_02}
Let $3/2 \le Y \le x$.  Let $g$ be a multiplicative function that for positive constants $\beta, c, c_1$ satisfies $|g(p)| \le \beta$,
\begin{equation}
\sum\limits_{w < p \le x} \left(|g(p)| - c\right)p^{-1} \ge -c_1, \quad Y \le w \le x, \notag
\end{equation}
on the primes.  Suppose, further, that the series
\begin{equation}
\sum\limits_q |g(q)| q^{-1} (\log q)^\gamma, \quad \gamma = 1 + c\beta(c+\beta)^{-1}, \notag
\end{equation}
taken over the prime-powers $q = p^k$ with $k \ge 2$, converges.

Then with
\begin{equation}
\lambda = \min\limits_{|t| \le T} \sum\limits_{Y < p \le x} \left(|g(p)| - \Re g(p)p^{it}\right)p^{-1}, \notag
\end{equation}
\begin{equation}
\sum\limits_{n \le x} g(n) \ll x(\log x)^{-1} \prod\limits_{p \le x} \left(1+|g(p)|p^{-1}\right) \left( \exp\left(-\lambda c(c+\beta)^{-1} \right) + T^{-1/2} \right) \notag
\end{equation}
uniformly for $Y$, $x$, $T > 0$, the implied constant depending at most upon $\beta, c, c_1$ and a bound for the sum of the series over higher prime-powers.
\end{theorem}   

\makeatletter{}

\noindent \emph{Remark}.  Since the partial sums $\sum_{p \le x} \left(|g(p)| - \Re g(p)p^{it} \right)p^{-1}$ are continuous in $t$ and nondecreasing in $x$, their divergence for each real $t$ ensures that
\begin{equation}
\sum\limits_{n \le x} g(n) = o\left(x(\log x)^{-1} \prod\limits_{p \le x} \left(1+|g(p)|p^{-1} \right) \right), \quad x \to \infty. \notag
\end{equation} 


As an example, Theorem \ref{JP_elliott_04_thm_02} applies to the renormalised Ramanujan function $|\tau(n)|n^{-11/2}$.

Since $\gamma > 1$, Lemmas \ref{elliott_duality_lem_02_02} and \ref{JP_elliott_04_lem_02} may be applied to the function $|g(n)|$.

Define
\begin{equation}
M(x) = \sum\limits_{n \le x} g(n), \quad N(x) = \sum\limits_{n \le x} g(n)\log n, \quad x \ge 2. \notag
\end{equation}
We switch between $M(x)$ and $N(x)$ as the harmonic analysis more readily lends itself.


\makeatletter{}

\begin{lemma} \label{JP_elliott_04_lem_04}
Assuming only that the multiplicative function $g$ is complex-valued,
\begin{equation}
\int_2^x \frac{|M(u)| \log u}{u^2} \, du \ll \int_2^x \frac{|N(u)|}{u^2} \, du \ll \int_1^x \frac{|M(u)| \log u}{u^2} \, du, \quad x \ge 2. \notag
\end{equation}
\end{lemma}     

\makeatletter{}

\noindent \emph{Proof}.  Integrating by parts:
\begin{equation}
M(u) = \frac{N(u)}{\log u} + \int_2^u \frac{N(w)}{w(\log w)^2} \, dw, \quad u \ge 2. \notag
\end{equation}
The first inequality in Lemma \ref{JP_elliott_04_lem_04} follows from the argument
\begin{align}
\int_2^x & \frac{\log u}{u^2} \int_2^u \frac{|N(w)|}{w(\log w)^2} \, dw \, du \notag \\
& = \int_2^x \frac{|N(w)|}{w(\log w)^2} \int_w^x \frac{\log u}{u^2} \, du\, dw \ll \int_2^x \frac{|N(w)|}{w^2 \log w} \, dw. \notag
\end{align}

Bearing in mind that $N(u)=0$ if $u<2$, a similar argument establishes the second inequality of Lemma \ref{JP_elliott_04_lem_04}. 


For ease of exposition, until otherwise stated, during the proof of Theorem \ref{JP_elliott_04_thm_02} $g$ will denote an exponentially multiplicative function, i.e., $g(p^k) = g(p)^k/k!$, $k \ge 2$, for which $|g(p)| \le \beta$.  The attendant Dirichlet series
\begin{equation}
G(s) = \sum\limits_{n=1}^\infty g(n)n^{-s} = \exp\left(\sum\limits_p g(p)p^{-s} \right), \quad G_0(s) = \sum\limits_{n=1}^\infty |g(n)|n^{-s}, \notag
\end{equation}
$s = \sigma + it$, absolutely convergent in the half-plane $\sigma = \Re(s)>1$, define there analytic functions for which we shall not require analytic continuation.

The multiplicativity of $g$ enables $x^{-1}N(x)$ and, therefore, $x^{-1}M(x)$, to be related to a weighted mean of itself, introducing a convenient smoothing.


\makeatletter{}

\begin{lemma} \label{JP_elliott_04_lem_05}
For each $K > 0$,
\begin{equation}
N(x) \ll x\int_1^x |M(u)| u^{-2} \, du + x \log\log x \sup\limits_{\sqrt{x} < u \le x} |M(u)| u^{-1} + x(\log x)^{-K} \notag
\end{equation}
uniformly for $x \ge 2$, the implied constant depending at most upon $\beta, K$.
\end{lemma}     

\makeatletter{}

\noindent \emph{Proof}.  For the duration of this proof only, set $y = x-x(\log x)^{-r}$, where $r$ is a positive number to be chosen presently.

Application of the Cauchy-Schwarz inequality and Lemma \ref{JP_elliott_04_lem_02} shows that
\begin{align}
|N(x) - N(w)| & \le (x-u + 1)^{1/2} \left(\sum\limits_{n \le x} \left(|g(n)|\log n\right)^2 \right)^{1/2} \notag \\
& \ll \left(x(\log x)^{-r}\right)^{1/2} \left(x(\log x)^{2+\beta^2} \right)^{1/2} \notag
\end{align}
uniformly for $y \le u \le x$.  Fixing $r$ at a sufficiently large value,
\begin{equation}
N(x) = (x-y)^{-1} \int_y^x N(u) \, du + O\left(x(\log x)^{-K} \right). \notag
\end{equation}

From the exponential multiplicativity of $g$,
\begin{equation}
N(u) = \sum\limits_{p \le u} g(p)\log p \, M(up^{-1}), \notag
\end{equation}
and with $z = (\log x)^{2r}$ the integral representing $N(x)$ has the corresponding representation
\begin{equation}
\sum\limits_{z < p \le x} g(p) p \log p \int_{y/p}^{x/p} M(v) \, dv + \sum\limits_{p \le z} g(p) \log p \int_y^x M(up^{-1}) \, du. \notag
\end{equation}

If $H$ denotes the supremum of $u^{-1}|M(u)|$ taken over the range $\sqrt{x} \le u \le x$, then the second sum is
\begin{equation}
\ll H\sum\limits_{p \le z} p^{-1} \log p \int_y^x u \, du \ll H x(x-y) \log z. \notag
\end{equation}

The first sum is
\begin{equation}
\ll \sum\limits_{z < p \le x} p \log p \int_{y/p}^{x/p} |M(v)| \, dv \ll \int_1^{x/z} |M(v)| \sum\limits_{y/v \le p \le x/v} p \log p \, dv. \notag
\end{equation}
For $v \le xz^{-1}$, $xv^{-1} \le \left(xv^{-1}z^{-1/2} \right)^2 \le \left((x-y)v^{-1} \right)^2$, so that by an old estimate of Hardy and Littlewood, or an application of Selberg's sieve, as in \cite{Elliott1979} Chapter 2,
\begin{equation}
\sum\limits_{y/v \le p \le x/v} p\log p \ll \frac{x}{v} \log \frac{x}{v} \left(\frac{x-y}{v}\right)\left(\log \left(\frac{x-y}{v}\right)\right)^{-1} \ll \frac{x(x-y)}{v^2}. \notag
\end{equation}
The corresponding contribution to $N(x)$ is
\begin{equation}
\ll x\int_1^x v^{-2} |M(v)| \, dv, \notag
\end{equation}
and Lemma \ref{JP_elliott_04_lem_05} is established. 

\makeatletter{}

\begin{lemma} \label{JP_elliott_04_lem_06}
\begin{equation}
G(s) \ll \prod\limits_{p \le x} \left(1+|g(p)|p^{-1} \right) e^{-\lambda} \left((\sigma - 1) \log x \right)^\beta \notag
\end{equation}
uniformly for $1+(\log x)^{-1} \le \sigma \le 2$.
\end{lemma}     

\makeatletter{}

\noindent \emph{Proof}.  For $2 \le w \le x$ define
\begin{equation}
\rho(w,t) = \sum\limits_{p \le w} \left(|g(p)| - \Re g(p)p^{-it} \right)p^{-1}. \notag
\end{equation}

From the Euler product representation
\begin{equation}
G(s) = \prod\limits_p \left(1+\sum\limits_{k=1}^\infty g(p^k)p^{-ks} \right) = \exp\left( \sum\limits_p g(p)p^{-s} \right) \notag
\end{equation}
we obtain the bounds
\begin{align}
G(s)G_0(\sigma)^{-1} & \ll \exp\left(-\sum\limits_p \left(|g(p)| - \Re g(p)p^{-it} \right)p^{-\sigma} \right) \notag \\
& \ll \exp(-\rho(y,t)), \quad y = \exp\left((\sigma - 1)^{-1} \right), \notag
\end{align}
in the strip $1 < \sigma \le 2$, since integration by parts and appeal to the Chebyshev inequality $\pi(x) \ll x/\log x$ shows the sums
\begin{equation}
\sum\limits_{p > y} p^{-\sigma}, \quad \sum\limits_{p \le y} \left(p^{-1} - p^{-\sigma} \right) \ll (\sigma-1) \sum\limits_{p \le y} p^{-1} \log p \notag
\end{equation}
to be uniformly bounded there.

For each $T$ the function
\begin{equation}
\sum\limits_{p \le w} 2|g(p)| p^{-1} - \min\limits_{|t| \le T} \rho(w,t) = \max\limits_{|t| \le T} \sum\limits_{p \le w} \left(|g(p)| + \Re g(p)p^{-it} \right)p^{-1} \notag
\end{equation}
is nondecreasing in $w$, hence
\begin{equation}
\min\limits_{|t| \le T} \rho(y,t) \ge \min\limits_{|t| \le T} \rho(x,t) - 2\sum\limits_{y < p \le x} |g(p)| p^{-1}, \notag
\end{equation}
uniformly for $1+(\log x)^{-1} \le \sigma \le 2$.  In particular,
\begin{equation}
G(s) \ll G_0(\sigma) \exp\left( 2\sum\limits_{y < p \le x} |g(p)| p^{-1} - \lambda \right) \notag
\end{equation}
over the same range of $\sigma$-values.

Apart from the factor $e^{-\lambda}$, this upper bound is
\begin{align}
& \ll \exp\left( \sum\limits_{p \le y} |g(p)| p^{-1} +  2\sum\limits_{y < p \le x} |g(p)| p^{-1} \right) \notag \\
& \ll \prod\limits_{p \le x} \left(1+|g(p)| p^{-1} \right) \exp\left(\beta \sum\limits_{y < p \le x} p^{-1} \right) \notag \\
& \ll \prod\limits_{p \le x} \left(1+|g(p)| p^{-1} \right) \left((\sigma-1) \log x \right)^\beta, \notag
\end{align}
completing the proof. 

\makeatletter{}

\begin{lemma} \label{JP_elliott_04_lem_07}
Under the hypothesis
\begin{equation}
\sum\limits_{w < p \le x} \left(|g(p)| - c\right)p^{-1} \ge -c_1, \quad Y \le w \le x, \notag
\end{equation}
the bound
\begin{equation}
G_0(\sigma) \ll \prod\limits_{p \le x} \left(1+|g(p)| p^{-1} \right) \left((\sigma-1)\log x\right)^{-c} \notag
\end{equation}
holds uniformly for $1+(\log x)^{-1} \le \sigma \le 1 + 2(\log Y)^{-1}$.

Moreover, with a suitable choice for $c_1$, the converse is valid.
\end{lemma}     

\makeatletter{}

\noindent \emph{Proof}.  With $\delta = 1+(\log x)^{-1}$, appeal to Euler product representations shows that
\begin{equation}
G_0(\delta)^{-1} G_0(\sigma) = \exp\left( -\sum\limits_p |g(p)| p^{-\delta} + \sum\limits_p |g(p)| p^{-\sigma} \right), \notag
\end{equation}
which by the Chebyshev bounds in Lemma \ref{JP_elliott_04_lem_06} is
\begin{equation}
\ll \exp\left( - \sum\limits_{w < p \le x} |g(p)| p^{-1} \right), \quad w = \exp\left((\sigma-1)^{-1} \right). \notag
\end{equation}

For $\sigma \le 1+(\log Y)^{-1}$ the asserted bound follows directly from the hypothesis and the estimate $\sum_{w < p \le x} p^{-1} = \log(\log x/\log w) + O(1)$.

For $1+(\log Y)^{-1} < \sigma \le 1+2(\log Y)^{-1}$, so that $Y^{1/2} \le w < Y$, we may replace $w$ in the second exponential by $Y$ and note that the sum $\sum_{Y^{1/2} < p \le Y} p^{-1}$ is bounded uniformly in $Y$.

The converse proposition may be obtained by reversing the steps. 

\makeatletter{}

\noindent \emph{Proof of Theorem \ref{JP_elliott_04_thm_02}}.  For $\alpha > 0$, $s = 1+\alpha + i\tau$, define
\begin{equation}
\theta(\alpha) = \left(\alpha \int_{1+\alpha-i\infty}^{1+\alpha + i\infty} \left| \frac{G'(s)}{s}\right|^2 \, d\tau \right)^{1/2} \ge 0. \notag
\end{equation}

An application of the Cauchy-Schwarz inequality shows that
\begin{align}
\int_2^x & \frac{|N(u)|}{u^2} \, du \le \left(\int_2^x \frac{|N(u)|^2}{u^3} \, du \int_2^x \frac{dv}{v} \right)^{1/2} \notag \\
& \ll \left( \log x \int_1^\infty |N(u)|^2 u^{-3-2\delta} \, du \right)^{1/2}, \quad \delta = 2(\log x)^{-1}. \notag
\end{align}
Since
\begin{equation}
s^{-1} G'(s) = \int_{-\infty}^\infty N(e^y)e^{-y\sigma} e^{-yi\tau} \, d\tau, \notag
\end{equation}
we may appeal to Plancherel's theorem
\begin{equation}
\int_1^\infty \frac{ |N(u)|^2}{u^{3+2\delta}} \, du = \frac{1}{2\pi} \int_{1+\delta -i\infty}^{1+\delta+i\infty} \left| \frac{G'(s)}{s} \right|^2 \, d\tau \notag
\end{equation}
from which, via Lemma \ref{JP_elliott_04_lem_04}
\begin{equation}
\int_{\sqrt{x}}^x \frac{|M(u)|}{u^2} \, du \ll \frac{1}{\log x} \int_2^x \frac{|M(u)| \log u}{u^2} \, du \ll \theta(\delta). \notag
\end{equation}

Thus
\begin{align}
\int_{3/2}^x \frac{|M(u)|}{u^2} \, du & \ll \int_{3/2}^x \frac{|M(u)|}{u^2} \int_u^{u^2} \frac{dy}{y \log y} \, du \notag \\
& \ll \int_{3/2}^{x^2} \frac{1}{y\log y} \int_{\sqrt{y}}^y \frac{|M(u)|}{u^2} \, du \, dy. \notag
\end{align}
If we assume that $g(p)=0$ for $p \le Y$, so that $M(u)=1$ over the range $u \le Y$, then this double integral is
\begin{equation}
\ll 1+\int_Y^{x^2} \theta\left(\frac{2}{\log y} \right) \frac{dy}{y\log y} \ll 1 + \int_{1/\log x}^{2/\log Y} \frac{\theta(\alpha)}{\alpha} \, d\alpha, \notag
\end{equation}
the last after the change of variable $y=e^{2/\alpha}$.

Appealing to Lemma \ref{JP_elliott_04_lem_05},
\begin{equation}
N(x) \ll x\int_{1/\log x}^{2/\log Y} \theta(\alpha) \alpha^{-1} \, d\alpha + x + x\log\log x \sup\limits_{\sqrt{x} < u \le x} |M(u)|u^{-1}. \notag
\end{equation}

In view of Lemmas \ref{elliott_duality_lem_02_02} and \ref{JP_elliott_04_lem_02}, with $P_x = \prod_{p \le x} \left(1+|g(p)| p^{-1} \right)$, the third of these bounding terms is $\ll x(\log x)^{-1} P_x \log\log x$.

The integral involving $\alpha$ we divide into two ranges, according to whether $\alpha \le \Delta(\log x)^{-1}$, where $1 \le \Delta \le 2\log x/\log Y$.

A further application of Plancherel's theorem together with the elementary bound $\sum_{p \le x} \log p \ll x$ shows that
\begin{align}
\int_{1+\alpha-i\infty}^{1+\alpha+i\infty} & \left| \frac{G'(s)}{sG(s)} \right|^2 \, d\tau = 2\pi \int_1^\infty \left| \sum\limits_{p \le u} g(p)\log p \right|^2 \frac{du}{u^{3+2\alpha}} \notag \\
& \ll \int_1^\infty u^{-1-2\alpha} \, du \ll \alpha^{-1}. \notag
\end{align}

In particular,
\begin{equation}
\int\limits_{|t| \le 1} \left| \frac{G'(s)}{G(s)} \right|^2 \, dt \ll \frac{1}{\alpha}, \quad \int\limits_{|t-m| \le 1} \left| \frac{G'(s)}{G(s)} \right|^2 \, dt \ll \frac{1}{\alpha}, \quad m \in \mathbb{Z}, \notag
\end{equation}
the last by application to the translated function with $g(p)p^{-im}$ in place of $g(p)$.

In view of the factorisation $G' = (G'/G)G$, for $T>0$
\begin{equation}
\int_{1+\alpha-i\infty}^{1+\alpha+i\infty} \left| \frac{G'(s)}{s} \right|^2 \, dt \ll \frac{1}{\alpha} \max\limits_{\substack{\sigma = \alpha \\ |t| \le T}} |G(s)|^2 + \frac{1}{\alpha} \sup\limits_{\substack{\sigma = \alpha \\ |t| > T}} |G(s)|^2 \sum\limits_{m > T} m^{-2}, \notag
\end{equation}
which, by applications of Lemmas \ref{JP_elliott_04_lem_06} and \ref{JP_elliott_04_lem_07} respectively, is
\begin{equation}
\ll \alpha^{-1} P_x^2 \left( e^{-\lambda}(\alpha \log x)^\beta + T^{-1/2} (\alpha \log x)^{-c} \right)^2. \notag
\end{equation}
The corresponding contribution to the integral over $\alpha$ is
\begin{align}
& \ll P_x e^{-\lambda} \int_{1/\log x}^{\Delta/\log x} (\log x)^\beta \alpha^{-1+\beta} \, d\alpha + T^{-1/2} P_x \int_{1/\log x}^{\Delta/\log x} (\log x)^{-c} \alpha^{-1-c} \, d\alpha \notag \\
& \ll P_x\left(e^{-\lambda}\Delta^\beta + T^{-1/2} \right). \notag
\end{align}

Emplying only Lemma \ref{JP_elliott_04_lem_07}, the contribution from the range $\Delta(\log x)^{-1} \le \alpha \le 2(\log Y)^{-1}$ is
\begin{equation}
\ll P_x \int_{\Delta/\log x}^{2/\log Y} (\log x)^{-c} \alpha^{-1-c} \, d\alpha \ll P_x \Delta^{-c}. \notag
\end{equation}

Before choosing $\Delta$ we note that we may assume $T \ge 1$, otherwise Theorem \ref{JP_elliott_04_thm_02} follows immediately from Lemmas \ref{elliott_duality_lem_02_02} and \ref{JP_elliott_04_lem_02}.  Since
\begin{align}
\Re & \left(   2^{-1} \int_{-1}^1 \sum\limits_{p \le x} \left(|g(p)| - g(p)p^{it} \right) p^{-1} \, dt \right) \notag \\
& = \sum\limits_{p \le x} |g(p)| p^{-1} + O\left(\sum\limits_{p \le x} (p \log p)^{-1} \right), \notag
\end{align}
$\lambda$ cannot exceed $\sum_{p \le x} |g(p)| p^{-1} + O(1)$, hence $\beta \log(\log x/\log Y) + c_2$ for an absolute constant, $c_2$.  The inequalities
\begin{equation}
e^{\lambda} \le e^{c_2} \left( \frac{\log x}{\log Y}\right)^\beta \le \left(\frac{2\log x}{\log Y} \right)^{c+\beta} \notag
\end{equation}
will be satisfied unless $Y$ exceeds a certain fixed power of $x$, depending upon $c$; in which circumstances Theorem \ref{JP_elliott_04_thm_02} again follows from Lemmas \ref{elliott_duality_lem_02_02} and \ref{JP_elliott_04_lem_02}.

We may thus choose $e^{-\lambda}\Delta^{\beta + c}=1$ to obtain the bound
\begin{equation}
N(x) \ll xP_x \left(\exp\left(-\lambda c(c+\beta)^{-1} \right) + T^{-1/2} + (\log x)^{-1} \log\log x \right) + x  \notag
\end{equation}
Moreover, in view of the above remark on the size of $\lambda$, we may omit the fourth error bound in favour of the first.

To widen the applicability of the estimate for $N(x)$ we note that by Lemma \ref{JP_elliott_04_lem_06}
\begin{equation}
\min\limits_{|t| \le T} \rho(u,t) \ge \min\limits_{|t| \le T} \rho(x,t) - 2\sum\limits_{u < p \le x} |g(p)| p^{-1}, \notag
\end{equation}
so that if $\sqrt{x} \le u \le x$, $Y \le \sqrt{x}$,
\begin{equation}
N(u) \ll uP_x\left( \exp\left(-\lambda c(c+\beta)^{-1} \right) + T^{-1/2} + (\log x)^{-1} \log\log x \right) \notag
\end{equation}
uniformly for $\sqrt{x} \le u \le x$.  Applied to the first integration by parts in Lemma \ref{JP_elliott_04_lem_04} this yields
\begin{equation}
M(x) \ll x(\log x)^{-1} P_x \left( \exp\left(-\lambda c(c+\beta)^{-1} \right) + T^{-1/2} + (\log x)^{-1} \log\log x \right), \notag
\end{equation}
since
\begin{equation}
\int_2^x \frac{N(w)}{w(\log w)^2} \, dw \ll P_x \int_2^{\sqrt{x}} \frac{dw}{\log w} + \sup\limits_{\sqrt{x} \le w \le x} \frac{|N(w)|}{w} \int_{\sqrt{x}}^x \frac{dw}{(\log w)^2}. \notag
\end{equation}
For $\sqrt{x} < Y \le x$ the improved bound on $M(x)$ also follows from  Lemmas \ref{elliott_duality_lem_02_02} and \ref{JP_elliott_04_lem_02}.

Taking into account the restrictions on $Y$ as necessary, we may now retrace our steps using the improved bound for $\sup_{\sqrt{x} \le u \le x} |M(u)|u^{-1}$ that the above argument implies and obtain the improved bound
\begin{equation}
N(x) \ll xP_x \left(\exp\left(-\lambda c(c+\beta)^{-1} \right) + T^{-1/2} + \left((\log x)^{-1} \log\log x \right)^2 \right). \notag
\end{equation}

Finitely many iterations of this procedure enable us to omit the resulting third term in the bound in favour of the first, and complete the proof of Theorem \ref{JP_elliott_04_thm_02} for uniformly bounded exponentially multiplicative functions that vanish on the primes not exceeding $Y$.

To remove the latter restriction, define exponentially multiplicative functions $f$, $h$ by
\begin{equation}
f(p) = \begin{cases} g(p) & \text{if} \ p \le Y, \\ 0 & \text{if} \ p > Y, \end{cases} \qquad h(p) = \begin{cases} 0 & \text{if} \ p \le Y, \\ g(p) & \text{if} \ p > Y,\end{cases} \notag
\end{equation}
so that $g = f*h$ and
\begin{equation}
\sum\limits_{n \le x} g(n) = \sum\limits_{u \le \sqrt{x}} f(u) \sum\limits_{v \le x/u} h(v) + \sum\limits_{v < \sqrt{x}} h(v) \sum\limits_{\sqrt{x} < u \le x/v} f(u). \notag
\end{equation}
The first of the doublesums we estimate by applying to $h$ the version of Theorem \ref{JP_elliott_04_thm_02} that we have obtained so far, the necessary uniformity in $u$ derived along the lines of that for $N(u)$, earlier.  The doublesum is then
\begin{equation}
\ll \sum\limits_{u \le \sqrt{x}} |f(u)| \frac{x}{u\log x} \prod\limits_{Y < p \le x/u} \left(1+\frac{|g(p)|}{p} \right)\left( \exp \left(-\lambda c(c+\beta)^{-1} \right) + T^{-1/2} \right) \notag
\end{equation}
which, since
\begin{equation}
\sum\limits_{u \le \sqrt{x}} |f(u)| u^{-1} \ll \prod\limits_{p \le Y} \left(1+|g(p)|p^{-1} \right), \notag
\end{equation}
falls within the bound of Theorem \ref{JP_elliott_04_thm_02}.

Set $\theta = 1/\log Y$.  Then, for $v < \sqrt{x}$,
\begin{equation}
\sum\limits_{\sqrt{x} < u \le x/v} |f(u)| \le \frac{1}{(\sqrt{x})^\theta} \sum\limits_{u \le x/v} |f(u)| u^\theta \ll \frac{x^{1-\theta/2}}{v\log x} \sum\limits_{u \le x} \frac{|f(u)|}{u^{1-\theta}}, \notag
\end{equation}
the second step by an application of Lemma \ref{elliott_duality_lem_02_02}; since if $q = p^k$, then $|f(q)| q^\theta \log q = \left(( |f(p)| p^\theta )^k \log p\right)/(k-1)! \le \left( (\beta Y^\theta )^k \log p \right) /(k-1)!$ and the function $M$ appearing in that Lemma is uniformly bounded in terms of $\beta$ alone.

The elementary inequality $e^t - 1 \le te^t$, valid for $t \ge 0$, shows that for primes $p$ not exceeding $Y$, $p^\theta \le 1+\theta p^\theta \log p \le 1 + e\theta \log p$; hence
\begin{align}
\sum\limits_{u \le x} |f(u)| u^{\theta-1} & \le \exp\left(\sum\limits_{p \le x} |f(p)| p^{\theta-1} \right) \notag \\
& \le \exp\left(\sum\limits_{p \le x} |f(p)| p^{-1} + \beta e \theta \sum\limits_{p \le Y} p^{-1} \log p \right) \notag \\
& \ll \prod\limits_{p \le Y} \left(1+ |g(p)| p^{-1} \right). \notag
\end{align}

The second of the above doublesums is
\begin{align}
& \ll \frac{x^{1-\theta/2}}{\log x} \sum\limits_{v < \sqrt{x}} \frac{|h(v)|}{v} \prod\limits_{p \le Y} \left(1+\frac{|g(p)|}{p} \right) \notag \\
& \ll \exp\left(\frac{-\log x}{2\log Y} \right) \frac{x}{\log x} \prod\limits_{p \le x} \left(1+ \frac{|g(p)|}{p} \right), \notag
\end{align}
an amount much smaller than the first term in the upper bound of Theorem \ref{JP_elliott_04_thm_02}.

To remove the restriction to exponential multiplicativity we express a general multiplicative function $g$ as a Dirichlet convolution $g_1 * g_2$, where $g_1$ is exponentially multiplicative and satisfies $g_1(p) = g(p)$.  Computation with Euler products shows that
\begin{align}
g_2(p^k) & = \sum\limits_{r=0}^k (r!)^{-1} (-g(p))^r g(p^{k-r}), \quad k \ge 2, \notag \\
g_2(p) & = 0. \notag
\end{align}

From the elementary inequality $a+b \le ab$, valid for reals satisfying $\min(a,b) \ge 2$,
\begin{equation}
4 \log n \le \prod\limits_{p^k \| n} (4 \log p^k). \notag
\end{equation}
For any nonnegative multiplicative function $h$, real $\gamma \ge 0$,
\begin{align}
\sum\limits_{n=1}^\infty h(n) (\log n)^\gamma & \le \prod\limits_p \left(1+ \sum\limits_{k=1}^\infty h(p^k) (4\log p^k)^\gamma \right) \notag \\
& \le \exp\left( 4^\gamma \sum\limits_q h(q)(\log q)^\gamma \right), \notag
\end{align}
where $q$ runs over the prime-powers, the final sum assumed to converge.

In the present circumstance
\begin{align}
\sum\limits_q q^{-1}  |g_2(q)| (\log q)^\gamma & \ll \notag \\
 \sum\limits_{p, t \ge 2} p^{-t} |g(p^t)| & (\log p)^t \left( 1+\sum\limits_{m=1}^\infty (m!p^m)^{-1} \beta^m (\log p^m)^\gamma \right) \notag \\
& + \sum\limits_{p,k \ge 2} \left( (k-1)! p^k \right)^{-1} \beta^k k^\gamma (\log p)^\gamma \notag \\
& \ll 1 + \sum\limits_{\substack{q = p^k \\ k \ge 2}} q^{-1} |g(q)| (\log q)^\gamma, \notag
\end{align}
so that the series $\sum_{n=1}^\infty n^{-1} |g_2(n)| (\log n)^\gamma$ with $\gamma = 1+c\beta(c+\beta)^{-1}$ converges.

The convolution decomposition of $g$ guarantees a representation
\begin{equation}
\sum\limits_{n \le x} g(n) = \sum\limits_{v \le x} g_2(v) \sum\limits_{u \le x/v} g_1(u). \notag
\end{equation}
Since $\sum_{\sqrt{x} < p \le x} p^{-1}$ is bounded uniformly in $x$, along earlier lines the contribution to the doublesum from the terms with $v \le \sqrt{x}$ is
\begin{equation}
\ll \sum\limits_{v \le \sqrt{x}} v^{-1} |g_2(v)| x(\log x)^{-1} P_x \left( \exp\left( -c\lambda (c+\beta)^{-1} \right) + T^{-1/2} \right), \notag
\end{equation}
within the asserted bound of Theorem \ref{JP_elliott_04_thm_02}.

The simple estimate
\begin{equation}
\sum\limits_{u \le x/v} |g_1(u)| \ll xv^{-1} P_x, \quad 1 \le v \le x, \notag
\end{equation}
shows  the terms with $\sqrt{x} < v \le x$ to contribute
\begin{equation}
\ll xP_x \sum\limits_{\sqrt{x} < v \le x} v^{-1} |g_2(v)| \ll xP_x (\log x)^{-1-\beta c(c+\beta)^{-1}}, \notag
\end{equation}
again within the bound of Theorem \ref{JP_elliott_04_thm_02}.

The proof of Theorem \ref{JP_elliott_04_thm_02} is complete.



\emph{Remarks}.  With an integration by parts the lower bound hypothesis on $|g(p)|$ in Theorem \ref{JP_elliott_04_thm_02} follows  directly from a uniform lower bound $\sum_{p \le y} p^{-1} |g(p)| \log p \ge c\log y - c_3$, $c_4 \le y \le x$.

We may braid $g$ with a Dirichlet character $\modd{D}$ without affecting the lower bound hypothesis of Theorem \ref{JP_elliott_04_thm_02} provided $D \le Y$.

\setcounter{example}{2}
\begin{example}\label{JP_elliott_04_exa_03}
Let $0 < c < 1$ and for $y \ge \exp\left(e^2\right)$ set $\beta(y) = (\log\log\log y)^{1/2}$.  Define a sequence $b_n$ of positive reals by $b_1 = 3/2$, $b_{n+1} = b_n\left(1+(\log b_n)^{-1}\right)$, $n \ge 1$.  The sequence is unbounded and for all sufficiently large values of $n$, $y = \left[b_n \left(\beta(b_{n+1}) - \beta(b_n)\right)\right]$ satisfies
\begin{equation}
b_n (2\log b_n)^{-4} \le y \le b_n(\log b_n)^{-2}. \notag
\end{equation}

Let $q_j$, $j = 1, 2, \dots$, enumerate the successive primes in the interval $(b_n, b_{n+1}]$.  Define $g(q_j)$ to be $-1$ if $1 \le j \le y$, 1 if $y < j \le c\left(\pi(b_{n+1}) - \pi(b_n)\right)$, and to be zero otherwise.  Note that $c\left(\pi(b_{n+1}) - \pi(b_n) \right) = \left(c+o(1)\right)b_n(\log b_n)^{-1}$.

Since
\begin{align}
\sum\limits_{j=1}^y \frac{g(q_j)}{q_j} & = \frac{y}{b_n} \left(1+O\left(\frac{1}{\log b_n}\right)\right)^{-1} \notag \\
& = \left(1+o(1)\right) \left(\beta(b_{n+1}) - \beta(b_n)\right), \notag
\end{align}
summing over the intervals with $b_n \le x$,
\begin{equation}
\sum\limits_{p \le x} \left(|g(p)| - g(p)\right) p^{-1} = \left(2+o(1)\right)(\log\log\log x)^{1/2}, \notag
\end{equation}
whilst
\begin{equation}
\sum\limits_{p \le x} |g(p)| \log p = \left(c+o(1)\right)x, \quad x \to \infty. \notag
\end{equation}

Suppose now that for some $|t|$ not exceeding $T$, a given fixed power of $\log x$,
\begin{equation}
\rho(x,t) = \sum\limits_{p \le x} \left(|g(p)| - \Re g(p)p^{-it} \right)p^{-1} < 2(\log\log\log x)^{1/2}. \notag
\end{equation}
Since
\begin{equation}
|g(p)| \left|1-\left(g(p)p^{it}\right)^2 \right|^2 \le 4 |g(p)| \left|1 - g(p)p^{it} \right|^2 \le 8\left(|g(p)| - \Re g(p)p^{it} \right), \notag
\end{equation}
\begin{equation}
\sum\limits_{p \le x} |g(p)| p^{-1} \left| 1-p^{2it} \right|^2 < 16(\log\log\log x)^{1/2}. \notag
\end{equation}
An application of I, Lemma 15 shows that $t \ll (\log x)^{-c/7}$.  With $w = \exp\left((\log x)^{c/7}\right)$,
\begin{equation}
\sum\limits_{p \le x} p^{-1} \left|p^{it}-1 \right| \ll \sum\limits_{p \le w} p^{-1} |t|\log p \ll 1. \notag
\end{equation}
In particular,
\begin{equation}
\rho(x,t) \ge \rho(w,t) = \left(2+o(1)\right)\beta(w) = \left(2+o(1)\right) \beta(x). \notag
\end{equation}

Hence
\begin{align}
\lambda & = \min\limits_{|t| \le T} \sum\limits_{p \le x} \left(|g(p)| - \Re g(p)p^{-it}\right) p^{-1} \notag \\
& = \left(2+o(1) \right)(\log\log\log x)^{1/2}, \quad x \to \infty, \notag
\end{align}
as required in the introduction to Theorem \ref{JP_elliott_04_thm_02}.
\end{example}

We include here S\"atze 1.1 and 1.2.2 of Wirsing \cite{wirsing1967}, to which we shall appeal in \S \ref{JP_elliott_04_sec_automorphic_forms_again}, and which may be compared with the present Theorem \ref{JP_elliott_04_thm_02}.

\makeatletter{}

\begin{lemma} \label{wirsing_1967_satz_01_01}
Let $\lambda(n)$ be a nonnegative real-valued multiplicative function that for a positive $\tau$ satisfies
\begin{equation}
\sum\limits_{p \le x} \frac{\log p}{p} \lambda(p) \sim \tau \log x, \quad x \to \infty. \notag
\end{equation}
Let the $\lambda(p)$ be uniformly bounded and the series $\sum \lambda(q)q^{-1}$, taken over the prime-powers $q = p^k$ with $k \ge 2$, converge.  Moreover, if $\tau \le 1$ then, with the same convention, let $\sum_{q \le x} \lambda(q) \ll x(\log x)^{-1}$ hold.

Then
\begin{equation}
\sum\limits_{n \le x} \lambda(n) \sim \frac{e^{-\kappa t}}{\Gamma(\tau)} \frac{x}{\log x} \prod\limits_{p \le x} \left(1 + \frac{\lambda(p)}{p} + \frac{\lambda(p^2)}{p^2} + \cdots \right), \quad x \to \infty, \notag
\end{equation}
where $\kappa$ denotes Euler's constant.
\end{lemma}  

\makeatletter{}

\begin{lemma} \label{wirsing_1967_satz_01_02_02}
Let $\lambda(n)$ satisfy the conditions of Lemma \ref{wirsing_1967_satz_01_01}, and let $g(n)$ be a real-valued multiplicative function that satisfies $|g(n)| \le \lambda(n)$ on the positive integers.

Then
\begin{equation}
\lim\limits_{x \to \infty} \sum\limits_{n \le x} g(n) \left( \sum\limits_{n \le x} \lambda(n) \right)^{-1} = \prod\limits_p \left( 1 + \sum\limits_{k=1}^\infty g(p^k)p^{-k} \right) \left(1+\sum\limits_{k=1}^\infty \lambda(p^k)p^{-k} \right)^{-1}, \notag
\end{equation}
where the product either converges properly to a nonzero limit, or improperly to zero.
\end{lemma} 

The asymptotic estimate for $\sum_{p \le x} p^{-1} \lambda(p)\log p$, rather than a lower bound, delivers an asymptotic estimate for the mean-value of $\lambda$.

It is not difficult to modify the proof of Theorem \ref{JP_elliott_04_thm_02}, including an analogue of Lemma \ref{elliott_duality_lem_02_02} for exponentially multiplicative functions, so that the theorem continues to hold under wider conditions.


\setcounter{theorem}{4}
\makeatletter{}

\begin{theorem} \label{JP_elliott_04_thm_05}
Theorem \ref{JP_elliott_04_thm_02} remains valid if, for some real $\delta$, $0 < \delta < 1$, the uniform bound $|g(p)| \le \beta$ on $g$ is replaced by the requirements
\begin{equation}
g(p) \ll p^{1/2 - \delta}, \qquad \sum\limits_{p \le x} |g(p)|^{1+\delta} p^{-1} \ll \log\log x, \notag
\end{equation}
\begin{equation}
\sum\limits_{x-y < p \le x} |g(p)| \log p \ll y \quad \text{uniformly for} \ x^{1-\delta} \le y \le x, \notag
\end{equation}
and
\begin{equation}
\sum\limits_{u < p \le v} |g(p)| p^{-1} \le \beta \log(\log v/\log u) + O(1) \notag
\end{equation}
uniformly for $v \ge u \ge 2$.
\end{theorem} 


In particular, Theorem \ref{JP_elliott_04_thm_05} may be applied to the Fourier coefficients of Maass forms.


\section{Automorphic Forms}\label{JP_elliott_04_sec_automorphic_forms_again}


We continue this study with three paradigmatic examples.

Let $a_n$, for integers $n$, be the coefficient function either of an elliptic holomorphic form of weight $k \ge 2$,
\begin{equation}
\sum\limits_{n=1}^\infty a_n e^{2\pi i nz}, \notag
\end{equation}
or of a Maass form
\begin{equation}
\sum\limits_{n \ne 0} a_n (2\pi y)^{1/2} K_w\left(2\pi |n| y\right) e^{2\pi i nx}, \notag
\end{equation}
i.e. a solution $f$, in an appropriate Hilbert space, of the Laplacian
\begin{equation}
-y^2\left( \frac{\partial^2 f}{\partial x^2} + \frac{\partial^2 f}{\partial y^2} \right) = \left( \frac{1}{4} - w^2 \right) f, \notag
\end{equation}
each a new cusp form attached to the action of a congruence subgroup $\Gamma_0(N)$ of $SL(2,\mathbb{Z})$ on the complex upper half-plane, eigenfunction of the appropriate Hecke operators and normalised to have $a_1 = 1$.

On the positive integers define the multiplicative function $g(n)$ to be 1 if $a_n \ne 0$, zero otherwise, and for each nonzero integer $D$ set
\begin{equation}
S_D(x) = \sum\limits_{\substack{n \le x \\ (n,D)=1}} g(n). \notag
\end{equation}


\makeatletter{}

\begin{theorem} \label{JP_elliott_04_thm_06}
Assume that for positive constants $c$, $c_0$, and all sufficiently large values of $x$
\begin{equation}
\sum\limits_{\substack{p \le x \\ a_p \ne 0}} \log p \ge cx \notag
\end{equation}
and
\begin{equation}
\sum\limits_{\substack{p \le x \\ a_p \ne 0}} p^{-1} \ge \frac{1}{2}\log\log x - c_0. \notag
\end{equation}

Then $S_D(x)$ lies between constant multiples of
\begin{equation}
x(\log x)^{-1} \prod\limits_{\substack{p \le x \\ (p,D)=1}} \left( 1+|g(p)| p^{-1} \right). \notag
\end{equation}
In particular $S_D(x) \gg x(\log x)^{-1/2}$, $x \ge 2$, the implied constants possibly depending upon $D$.

For each pair of mutually prime integers $a$, $D$,
\newlength{\Myll}
\settowidth{\Myll}{$S_D(x)^{-1}$}
\begin{equation}
\gamma(a) = \lim\limits_{x \to \infty} \mathop{ S_D(x)^{-1} \sum \makebox[\Myll][l]{$g(n)$}} \limits_{\substack{n \le x \\ n \equiv a \modd{D}}}   \notag
\end{equation}
exists.  Moreover,
\begin{enumerate}
\item $\gamma(a) = \varphi(D)^{-1}$ uniformly in $(a,D)=1$ unless there is a quadratic character $\chi \modd{D}$ for which the series $\sum g(p)p^{-1}$, taken over the primes with $\chi(p)=-1$, converges.
\item In that exceptional case
\begin{equation}
\gamma(a) = \varphi(D)^{-1} \left( 1 + \chi(a)\prod\limits_{\substack{g(p)=1 \\ \chi(p) = -1}} \psi_p \prod\limits_{\substack{g(p)=0 \\ \chi(p)=-1}} \psi_p \right) \notag
\end{equation}
with
\begin{equation}
\psi_p = \left( 1 + \sum\limits_{k=1}^\infty g(p^k)\chi(p^k)p^{-k} \right)\left( 1 + \sum\limits_{k=1}^\infty g(p^k)p^{-k} \right)^{-1}, \notag
\end{equation}
and is nonzero unless, for the primes not dividing $D$, $g(p)=1$ if and only if $\chi(p) = +1$.
\item In that further exceptional case $\gamma(a) = 2\varphi(D)^{-1}$ when $\chi(a)=1$, and is zero otherwise.
\end{enumerate}

In the cases $(2)$ and $(3)$
\begin{equation}
S_D(x) = \left(1+o(1)\right) Ax(\log x)^{-1/2}, \quad x \to \infty, \notag
\end{equation}
for a certain positive constant $A$.
\end{theorem}


\emph{Remarks}.  For holomorphic forms Theorem \ref{JP_elliott_04_thm_06} is in accord with the result of Serre considered in \S \ref{JP_elliott_04_sec_automorphic_forms}, the existence of an exceptional quadratic character corresponding to the form having complex multiplication.  The required conditions are all satisfied.

For Maass forms the argument in \S \ref{JP_elliott_04_sec_automorphic_forms} shows that, in the notation of Kim and Shahidi loc. cit., provided the underlying representation of $GL_2(A_\mathbb{Q})$ is not of dihedral, tetrahedral or octahedral type, the conditions are also satisfied.


\makeatletter{}

\noindent \emph{Proof of Theorem \ref{JP_elliott_04_thm_06}}.  The bounds on $S_D(x)$ follow from the first three lemmas in \S \ref{JP_elliott_04_sec_multiplicative_functions}.  It is only here that we appeal to the assumption that $\sum_{p \le x} \log p$, taken over the primes for which $a_p$ does not vanish, is at least $cx$ in size.  For the remainder of the present proof it will suffice that $\sum_{p \le x} p^{-1} \log p \ge c\log x$, with the same constraint on the primes, be satisfied.

From the representation of $g$, on arithmetic progressions, in terms of Dirichlet characters and by the remark following Theorem \ref{JP_elliott_04_thm_02}, $\gamma(a)$ will exist and have the value $\varphi(D)^{-1}$ unless, for some nonprincipal character $\chi \modd{D}$ and real $t$, the series
\begin{equation}
\sum\limits_p g(p)  \left(1- \Re \chi(p) p^{it} \right) p^{-1} \notag
\end{equation}
converges.

If the character has order $r$, appeal to the inequality $\left| 1-z^r \right| \le r|1-z|$, $|z| \le 1$, shows the series
\begin{equation}
\sum\limits_p g(p) \left( 1-\Re p^{irt} \right) p^{-1} \notag
\end{equation}
to converge and, after I, Lemma 15, that $t=0$.

Moreover, by I, Lemma 3,
\begin{align}
\sum\limits_{p \le x} g(p)p^{-1} & \le \Re \sum\limits_{j=1}^r r^{-1} \sum\limits_{p \le x} \chi(p)^j p^{-1} + O(1) \notag \\
& \le \left(r^{-1} + o(1) \right) \log\log x, \quad x \to \infty, \notag
\end{align}
contradicting an initial assumption unless $r=2$.

Likewise, there can be at most one exceptional quadratic character, for an inequivalent second such, $\chi_1$, would guarantee
\begin{align}
\sum\limits_{p \le x} g(p) p^{-1} & \le \sum\limits_{p \le x} 4^{-1} \left( 1 + \chi(p) \right) \left( 1 + \chi_1(p) \right) p^{-1} + O(1) \notag \\
& \le \left(1/4 + o(1) \right) \log\log x, \quad x \to \infty, \notag
\end{align}
again a contradiction.

In view of the identity
\begin{equation}
\mathop{ \phantom{p^{-1}} \sum p^{-1} }\limits_{\substack{p \le x \\ g(p)=0, \, \chi(p)=1}} = \sum\limits_{\substack{p \le x \\ \chi(p)=1}} p^{-1} - \left( \sum\limits_{\substack{p \le x \\ g(p)=1}} p^{-1} - \mathop{ \phantom{p^{-1}} \sum p^{-1} }\limits_{\substack{p \le x \\ g(p)=1, \, \chi(p)=-1}}  \right) \notag
\end{equation}
and the lower bound hypothesis on $\sum_{p \le x} g(p)p^{-1}$, there is a set of primes $q$ for which $\sum q^{-1}$ converges and outside of which $g(p)=1$ if and only if $\chi(p)=1$.

In particular,
\begin{equation}
\sum\limits_{q \le x} q^{-1} \log q \le \sum\limits_{q \le \log x} q^{-1} \log q + (\log x) \sum\limits_{q > \log x} q^{-1} = o(\log x), \notag
\end{equation}
hence
\begin{equation}
\sum\limits_{p \le x} g(p)p^{-1} \log p = \sum\limits_{\substack{p \le x \\ \chi(p)=1}} p^{-1} \log p + o(\log x) = \left(1/2 + o(1) \right) \log x, \quad x \to \infty. \notag
\end{equation}

Lemmas \ref{wirsing_1967_satz_01_01} and \ref{wirsing_1967_satz_01_02_02} deliver the asymptotic estimates
\begin{equation}
S_D(x) = \left(1+o(1) \right) \frac{e^{-\kappa/2}}{\Gamma(1/2)} \frac{x}{\log x} \prod\limits_{\substack{p \le x \\ (p,D)=1}} \left( 1+ \sum\limits_{m=1}^\infty \frac{g(p^m)}{p^m} \right), \notag
\end{equation}
\begin{equation}
\lim\limits_{x \to \infty} S_D(x)^{-1} \sum\limits_{n \le x} g(n)\chi(n) = \prod\limits_p \psi_p, \notag
\end{equation}
respectively.  Here $\kappa$ denotes Euler's constant and our essential characterisation of the primes on which $g(p)=1$ allows us to deduce that
\begin{equation}
S_D(x) = \left(1+o(1) \right) Ax(\log x)^{-1/2}, \quad x \to \infty. \notag
\end{equation}

The remaining assertions of Theorem \ref{JP_elliott_04_thm_06} follow rapidly.

This completes the proof of Theorem \ref{JP_elliott_04_thm_06}.


Let $a_n$ be the Fourier coefficient function of a Maass cusp form corresponding to the action of $SL(2, \mathbb{Z})$ upon the complex upper half-plane, eigenfunction of the appropriate Hecke operators, normalised to have $a_1 = 1$.


\makeatletter{}

\begin{theorem} \label{JP_elliott_04_thm_07}
Let $S(x)$ denote the number of integers $n$, not exceeding $x$, for which $a_n$ does not vanish.  Then $S(x) \gg x(\log x)^{-1/2}$ and, with $24000\gamma=1$,
\begin{equation}
S(x)^{-1} \sum\limits_{\substack{ n \le x \\ a_n < 0}} 1 = 1/2 + O\left( (\log x)^{-\gamma} \right), \quad x \ge 2. \notag
\end{equation}
There is a like result for the $a_n$ that are positive.
\end{theorem} 


\emph{Remarks}.  $S(x)$ coincides with the function $S_1(x)$ of Theorem \ref{JP_elliott_04_thm_06}.

The positive and negative $a_n$ are uniformly distributed amongst those that are nonzero.

Provided the coefficients are real, and we adjust the value of $\gamma$ as necessary, various versions of Theorem \ref{JP_elliott_04_thm_07} may be obtained for holomorphic as well as Maass forms, moreover $SL(2,\mathbb{Z})$ may be replaced by a congruence subgroup $\Gamma_0(N)$.  In particular, Theorem \ref{JP_elliott_04_thm_07} holds for the values of Ramanujan's $\tau$-function.

We may restrict summations to be arithmetic progressions, with the same interests as those attached to Theorem \ref{JP_elliott_04_thm_06}.

For the duration of the proof of the present theorem only, we shall denote $\log\log x$ by $L$.

We first furnish a supply of negative coefficients.


\makeatletter{}

\begin{lemma} \label{JP_elliott_04_lem_10}
\begin{equation}
\sum\limits_{\substack{p \le x \\ a_p < 0}} p^{-1} \ge 6^{-3}L + O(1), \quad x \ge e. \notag
\end{equation}
\end{lemma} 

\makeatletter{}

\noindent \emph{Proof of Lemma \ref{JP_elliott_04_lem_10}}.  We recall the estimates
\begin{equation}
\sum\limits_{p \le x} |a_p|^2 p^{-1} = L + O(1), \quad \sum\limits_{p \le x} |a_p|^4 p^{-1} = 2L + O(1), \quad x \ge e, \notag
\end{equation}
from the introductory section.

For any real $\lambda > 0$,
\begin{align}
\sum\limits_{\substack{p \le x \\ |a_p| \le \lambda}} |a_p|^2 p^{-1} & \ge L - \lambda^{-2} \sum\limits_{p \le x} |a_p|^4 p^{-1} + O(1) \notag \\
& \ge \left( 1 - 2\lambda^{-2} \right)L + O(1). \notag
\end{align}
Hence
\begin{equation}
\sum\limits_{\substack{ p \le x \\ |a_p| \le \lambda}} |a_p| p^{-1} \ge \lambda^{-1} \left( 1 - 2\lambda^{-2} \right) L + O(1). \notag
\end{equation}
The coefficient of $L$ is maximized by $\lambda^2 = 6$.  In particular,
\begin{equation}
\sum\limits_{p \le x} |a_p| p^{-1} \ge 2\bigl(3\sqrt{6}\bigr)^{-1}L + O(1), \quad x \ge e. \notag
\end{equation}

The function $f(s) = \sum_{n=1}^\infty a_n n^{-s}$ has a holomorphic continuation over the whole plane.  We need only that it is holomorphic in a proper disc $|s-1| < c$.  Assuming it to have a zero of order $r$ at $s=1$, consideration of the function $\log\left( f(s)\zeta(s)^r \right)$ and its derivative shows that
\begin{equation}
\sum\limits_{p \le x} a_p \, p^{-1} = -rL + O(1), \quad x \ge e. \notag
\end{equation}
As a consequence
\begin{equation}
2\sum\limits_{\substack{p \le x \\ a_p < 0}} |a_p| p^{-1} = \sum\limits_{p \le x} \left(|a_p| - a_p \right)p^{-1} \ge \sum\limits_{p \le x} |a_p| p^{-1} + O(1). \notag
\end{equation}

For any $\nu > 0$,
\begin{align}
\sum\limits_{\substack{p \le x \\ a_p < 0, \, |a_p| \le \nu}} |a_p| p^{-1} & \ge \bigl(3\sqrt{6}\bigr)^{-1} L - \nu^{-1} \sum\limits_{p \le x} |a_p|^2 p^{-1} + O(1) \notag \\
& = \left( \bigl(3\sqrt{6} \bigr)^{-1} - \nu^{-1} \right)L + O(1). \notag
\end{align}
In particular,
\begin{equation}
\sum\limits_{\substack{p \le x \\ a_p < 0}} p^{-1} \ge \nu^{-1} \left( \bigl(3\sqrt{6}\bigr)^{-1} - \nu^{-1} \right)L + O(1), \quad x \ge e. \notag
\end{equation}
The choice $\nu = 6\sqrt{6}$ completes the proof of Lemma \ref{JP_elliott_04_lem_10}.



\makeatletter{}

\noindent \emph{Proof of Theorem \ref{JP_elliott_04_thm_07}}.  Define the multiplicative function $g$ by
\begin{equation}
g(p^k) = \begin{cases} 1 & \text{if} \ a_{p^k} > 0, \\ -1 & \text{if} \ a_{p^k} < 0, \\ 0 & \text{if} \ a_{p^k}=0. \end{cases} \notag
\end{equation}
Then $g(n)$ has the same sign as $a_n$ and it will suffice to prove that restricted to the integers for which $a_n$ does not vanish, $g$ has mean-value zero.

With $c = 1/2$, $\beta = 1$, $T  = (\log x)^4$,
\begin{equation}
\rho = \min\limits_{|t| \le T} \sum\limits_{p \le x} \left( |g(p)| - \Re g(p)p^{it} \right) p^{-1}, \notag
\end{equation}
Theorem \ref{JP_elliott_04_thm_02} delivers the bound
\begin{equation}
S(x)^{-1} \sum\limits_{n \le x} g(n) \ll e^{-\rho/3}, \quad x \ge e. \notag
\end{equation}

Let $t$ be a value for which the minimum defining $\rho$ is attained.  If, for some positive $\delta$ not exceeding $1/8$, $\rho \le \delta L$, then
\begin{equation}
\sum\limits_{p \le x} |g(p)| p^{-1} \left| 1 - \left( g(p) p^{it} \right)^2 \right|^2 \le 4 \sum\limits_{p \le x} |g(p)| p^{-1} \left| 1- p^{2it} \right|^2 \notag
\end{equation}
in turn $\le 8 \rho \le 8 \delta L$.  Application of I, Lemma 15 shows that either $2|t| \le (\log x)^{-1}$ or
\begin{equation}
\sum\limits_{p \le x} |g(p)| p^{-1} \le 4( 8\delta)^{1/3} -6\log |t| + O(1), \notag
\end{equation}
the implied constant depending at most upon $\delta$.  Recalling from the introduction the lower bound estimate
\begin{equation}
\sum\limits_{\substack{ p \le x \\ a_p \ne 0}} p^{-1} \ge 2^{-1} L + O(1), \notag
\end{equation}
we see that provided $4(8\delta)^{1/3} < 1/2$,
\begin{equation}
|t| (\log x)^\alpha \quad \text{with} \quad \alpha = \frac{1}{6} \left(\frac{1}{2} - 4(8\delta)^{1/3} \right) \notag
\end{equation}
is bounded uniformly in $x \ge e$.

Set $Z = \exp\left( (\log x)^\alpha \right)$.  Then
\begin{equation}
\sum\limits_{p \le Z}  \left| 1-p^{it} \right| p^{-1} \ll |t| \sum\limits_{p \le Z} p^{-1} \log p \ll |t| \log Z \ll 1, \notag
\end{equation}
and
\begin{align}
\rho & \ge \sum\limits_{p \le Z} \left( |g(p)| - \Re g(p) p^{it} \right) p^{-1} \ge \sum\limits_{p \le Z} \left( |g(p)| - g(p) \right) p^{-1} + O(1) \notag \\
& \ge 2 \cdot 6^{-3} \log\log Z + O(1) = 2\cdot 6^{-3} \alpha L + O(1). \notag
\end{align}

We may therefore choose for the role of $3\gamma$ the minimum of $\delta$ and $2 \cdot 6^{-3} \alpha$.

With $8 \delta = 10^{-3}$, $2\cdot 6^{-3}\alpha = 6^{-4} \cdot 5^{-1} > 10^{-3} \cdot 8^{-1} = \delta$.

The proof of Theorem \ref{JP_elliott_04_thm_07} is complete.


Between the function $g$ in Theorem \ref{JP_elliott_04_thm_06} and the coefficient function $n \mapsto |a_n|^2$ lies the arithmetic function $n \mapsto |a_n|$.  Even if the corresponding Dirichlet series $\sum_{n=1}^\infty |a_n| n^{-s}$ were to be analytically continuable, an essential singularity would be expected at the point $s = (k+1)/2$, where $k$ is the weight of the form.


\makeatletter{}

\begin{theorem} \label{JP_elliott_04_thm_08}
To each Maass cusp form satisfying the conditions of Theorem \ref{JP_elliott_04_thm_07} there corresponds an integer $c$ with the property that for any pair of mutually prime positive integers $a$, $D$,
\begin{equation}
\lim\limits_{x \to \infty} \left( \sum\limits_{n \le x} |a_n| \right)^{-1} \sum\limits_{\substack{n \le x \\ n \equiv a \modd{D}}} |a_n| = \varphi(D)^{-1} \notag
\end{equation}
holds unless there is a nonprincipal character $\chi \modd{D}$, of order at most $c$, for which the series
\begin{equation}
\sum\limits_{\substack{a_p \ne 0 \\ \chi(p) \ne 1}} |a_p| p^{-1}, \notag
\end{equation}
taken over the primes, converges.
\end{theorem} 

\emph{Remark}.  As for Theorem \ref{JP_elliott_04_thm_07}, the theorem remains valid if we replace $SL(2,\mathbb{Z})$ by $\Gamma_0(N)$ and, with appropriate renormalisation, allow the form to be holomorphic.

\makeatletter{}

\noindent \emph{Proof of Theorem \ref{JP_elliott_04_thm_08}}.  Since an analogous result is implicit in Theorem \ref{JP_elliott_04_thm_06}, we confine ourselves to the details of applying Theorem \ref{JP_elliott_04_thm_05} rather than Theorem \ref{JP_elliott_04_thm_02}.

With $g(p) = |a_p|$ the only condition that does not follow rapidly from the remarks in \S \ref{JP_elliott_04_sec_automorphic_forms} is the requirement that for some $\delta$, $0 < \delta < 1$, $\sum_{x-y < p \le x} |g(p)| \log p \ll y$ uniformly for $x^{1-\delta} \le y \le x$.

The analytic continuation over the simple pole at $s=1$ and into the half-plane $\Re(s) > 1/2$ of the Rankin-Selberg product $L$-function $\sum_{n=1}^\infty |a_n|^2 n^{-s}$ guarantees a positive $\delta$, $0 < \delta < 1/2$, for which
\begin{equation}
\sum\limits_{n \le x} |a_n|^2 = Ax + O\left(x^{1-2\delta} \right). \notag
\end{equation}
By subtraction
\begin{equation}
\sum\limits_{x-y < n \le x} |a_n|^2 = Ay + O\left( y^{1-\theta} \right), \quad \theta = \delta(1-\delta)^{-1}, \notag
\end{equation}
uniformly for $x^{1-\delta} \le y \le x$.  This is sufficient that we may apply the version of Selberg's sieve given in Chapter 2 of the first author's volume \cite{Elliott1979}, and obtain the bound
\begin{equation}
\sum\limits_{x-y < p \le x} |a_p|^2 \ll y(\log y)^{-1} \notag
\end{equation}
uniformly for $y \ge x^{1-\delta}$, $x \ge 2$.  That $a_n$ is a multiplicative function of $n$ is vital.  A similar application may be found in Elliott \cite{elliott2012CLTeigenforms}.

Appeal to the Cauchy-Schwarz inequality yields
\begin{equation}
\sum\limits_{x-y < p \le x} |a_p| \ll y(\log y)^{-1} \notag
\end{equation}
with the same uniformity, and the required condition is satisfied.

In order to apply Lemmas 12 and 15 from the taxonomy section of paper I we require that the function $h(p)$ there, here played by $a_p$, not exceed 1.  If, for some $\beta > 0$, and with $L = \log\log x$, $\sum_{p \le x} |a_p| p^{-1} \left| 1-\chi(p)p^{it} \right|^2 \le \beta L$, then we restrict the primes to those for which $|a_p|$ does not exceed $\sqrt{6}$.  It was shown during the proof of the present Lemma \ref{JP_elliott_04_lem_10} that over such primes $\sum_{p \le x} |a_p| p^{-1}$ is at least $2\bigl(3\sqrt{6}\bigr)^{-1} L + O(1)$.  We may proceed with $\beta$ replaced by $\beta/\sqrt{6}$, $h(p) = |a_p| / \sqrt{6}$ unless $|a_p| > \sqrt{6}$, when we set $h(p) = 0$.

This completes the outline proof of Theorem \ref{JP_elliott_04_thm_08}. 


\emph{Uniformities}.  In general, the argument employed in \S \ref{JP_elliott_04_sec_automorphic_forms} to guarantee nonvanishing Fourier coefficients $a_p$ may be expected to encounter poles of order greater than $2$ and deliver
\begin{equation}
\sum\limits_{\substack{p \le x \\ a_p \ne 0}} p^{-1} \ge \alpha \log\log x + O(1) \notag
\end{equation}
for a positive $\alpha$ possibly smaller than $1/2$.  The exceptional characters in Theorem \ref{JP_elliott_04_thm_06} may then have orders up to $\left[ \alpha^{-1} \right]$, and be better classified in terms of exceptional representations of the group $GL_2(A_\mathbb{Q})$, and so on.

With a view towards applications it is natural to consider Fourier coefficients of automorphic forms on finite arithmetic progressions and allow the differences to vary.  Within the aesthetic of the second section this has serious repercussions since, as was shown in the final example of paper I, with $g$ the M\"obius function results bear upon the distribution of primes.

Assume that the function $g$ in Theorem \ref{JP_elliott_04_thm_06} satisfies the bound
\begin{equation}
\sum\limits_{D<p \le y} g(p)p^{-1} \ge 1/2 \sum\limits_{D < p \le y} p^{-1} + O(1) \notag
\end{equation}
that is derived in \S \ref{JP_elliott_04_sec_automorphic_forms}.  Given $x \ge D$, $\gamma > 0$, confining interest to the interval $(D,y]$ and without assuming $D$ to be fixed, I, Theorem 3 delivers the estimate
\begin{align}
 \sum\limits_{\substack{n \le y \\ n \equiv a \modd{D}}} g(n)  = \frac{1}{\varphi(D)} & \sum\limits_{\substack{n \le y \\ (n,D)=1}} g(n) + \sum\limits_j \frac{\overline{\chi_j}(a)}{\varphi(D)} \sum\limits_{n \le y} g(n) \chi_j(n) \notag \\
 & + O\Biggl( \frac{y}{\varphi(D)\log y} \prod\limits_{\substack{p \le y \\ (p,D)=1}} \left(1+\frac{g(p)}{p} \right)\left( \frac{\log D}{\log y} \right)^\eta \Biggr) \notag
\end{align}
with $(c+1)\eta = 10^{-5} c^6$, the number of characters bounded in terms of $c$ alone, each character of order at most $20c^{-1}$, uniformly for $x^\gamma \le y \le x$.  Moreover, with $Z = \exp\left(\log D(\log x/\log D)^{1/30} \right)$ and the same uniformity, each exceptional sum satisfies the bound
\begin{equation}
\sum\limits_{n \le y} g(n) \chi_j(n) \ll \frac{y}{\log y} \prod\limits_{\substack{p \le y \\ (p,D)=1}} \left(1+ \frac{g(p)}{p} \right) \exp\left( -\frac{c}{c+1}  \sum\limits_{D < p \le Z} \frac{g(p)\left(1-\Re \chi_j(p) \right)}{p} \right). \notag
\end{equation}
Arguing much as in the proof of Theorem \ref{JP_elliott_04_thm_06}, we may reduce the set of characters $\chi_j$ until it contains at most a single real character.

However, there is no immediate reason that $g(p) = \chi_1(p)=1$ might not hold on successive stretches of primes with the defining moduli of the various real characters $\chi_1$ increasing.  This raises an interesting question:

Given a constant $B$ and a (possibly general) character $\chi$ for which
\begin{equation}
\sum\limits_{\substack{ p \le x \\ a_p \ne 0}} \left( 1 - \Re \chi(p) \right) p^{-1} \le B, \notag
\end{equation}
is there a value of $x$ that will force the background representation of the group $GL_2(A_\mathbb{Q})$ or $GL_2(A_F)$, as the case may be, to be of an exceptional type; and if so, how large does $x$ need to be as a function of the parameters of the character, in our case the modulus $D$?

Note that we do not ask to identify the form, only the class to which it belongs.

To complete this study and to illustrate a uniformity with respect to the modulus $D$, we overfly the apparatus.

\begin{example}\label{JP_elliott_04_exa_04}
Let $a_n$ be the Fourier coefficient function of a Maass cusp form attached to the action of $SL(2,\mathbb{Z})$ on the complex upper half-plane, eigenfunction of the appropriate Hecke operators.  Then for each pair of mutually prime positive integers $b$, $D$,
\begin{equation}
\mathop{\phantom{1} \sum 1}\limits_{\substack{n \le x, \,a_n \ne 0 \\ n \equiv b \modd{D}}} = \frac{1}{\varphi(D)} \left(1 + O\left( \left(\frac{\log D}{\log x} \right)^{1/49} \right) \right) \mathop{\phantom{1} \sum 1}\limits_{\substack{n \le x, \, a_n \ne 0 \\ (n,D)=1}}, \quad x \to \infty, \notag
\end{equation}
uniformly for $(b,D)=1$, $1 \le D \le x$.

\noindent \emph{Proof}.  We may assume that $D \ge 2$.  A rapid argument may appeal to I, Theorem 3.  A better error term applies I, Theorem 1 directly.  In the notation of the proof of Theorem \ref{JP_elliott_04_thm_06}, this reduces the problem to the estimation of a uniformly bounded number of sums $\sum_{n \le x} g(n)\chi_j(n)$, where $\chi_j$ is a nonprincipal Dirichlet character $\modd{D}$.

In turn, application of Theorem \ref{JP_elliott_04_thm_02} with $Y=D$ further reduces the problem to the provision of a lower bound for the sum
\begin{equation}
\rho = \sum\limits_{D < p \le x} g(p) \left( 1-\Re \chi(p) p^{it} \right)p^{-1}, \notag
\end{equation}
where $t$ is a real value in the interval $|t| \le \max\left(D, (\log x/\log D)^4 \right)$.

Rather than appeal to the taxonomy section of paper I, it is worthwhile to establish an analogue of I, Lemma 3.

Denote by $f(s)$ the function $\sum_{n=1}^\infty a_n^2 \chi(n) n^{-s}$, the Rankin-Selberg product of the $L$-function attached to the Maass form and the same $L$-function braided with a nonprincipal Dirichlet character $\modd{D}$.

Assume first that $\chi$ is primitive.  It follows from the work of Shimura \cite{shimura1975holomorphy}, c.f. Goldfeld \cite{Goldfeld2006}, that $f(s)$ has an analytic continuation over the whole complex plane and, after an appeal to its functional equation together with an application of the Phragm\'en-Lindel\"of theorem, c.f. Elliott \cite{elliott2012centrallimitramanujantau}, satisfies $f(s) \ll \left( D(2 + |s| )\right)^{c_0}$ for some absolute constant $c_0$, uniformly for $\Re(s) \ge 3/4$.

Examination of its Euler product shows that $f$ continues to satisfy such a bound whether $\chi$ is primitive or not.

The Euler product representation of $f$ factorises to
\begin{align}
\prod\limits_{p \mid D} \left( 1- p^{-2s} \right)^{-1} & \zeta(2s) f(s) \notag \\
& = \prod\limits_{(p,D)=1} \left( 1 - \alpha_p^2 \chi(p) p^{-s} \right)^{-1} \left( 1 - \chi(p) p^{-s} \right)^{-2} \left( 1 - \beta_p^2 \chi(p) p^{-s} \right)^{-1} \notag
\end{align}
where, after Kim and Shahidi loc. cit., the components of the local Satake parameter $\text{diag}(\alpha_p, \beta_p)$ satisfy $|\alpha_p| \le p^{1/9}$, $|\beta_p| \le p^{1/9}$.  Since $\alpha_p \beta_p = 1$, $\alpha_p + \beta_p = a_p$, taking logarithms and reassembling,
\begin{equation}
f(s)^{-1} = \exp\left( - \sum\limits_{(p,D)=1} a_p^2 \chi(p) p^{-s} + O(1) \right) \notag
\end{equation}
uniformly for $\Re(s) \ge 3/4$, and
\begin{equation}
\sigma - 1 \ll |f(s)| \ll (\sigma-1)^{-1} \notag
\end{equation}
uniformly for $1 < \sigma = \Re(s) \le 2$.

For a suitably chosen point $s_0$ in the strip $1 < \Re(s) < 2$, this enables the application of a Borel-Carath\'eodory inequality to $f(s) f(s_0)^{-1}$ in the style of Landau, as in Elliott \cite{elliottMFoAP6}, Lemma 14, to obtain the upper bound
\begin{equation}
\left| \frac{f(\sigma_1 + it)}{f(\sigma_2 + it)} \right| \ll \left( D(2+|s|) \right)^{c_1(\sigma_2 - \sigma_1)} \notag
\end{equation}
with a constant $c_1$ independent of $D$, uniformly for $1 < \sigma_1 \le \sigma_2$.

We may now adapt the argument of I, Lemma 3 to show that
\begin{equation}
\Re \sum\limits_{y < p \le w} a_p^2 \chi(p) p^{-1-it} \le \log(\log T/\log D) + c_2 \notag
\end{equation}
uniformly for $\Re(s) \ge 1$, $|t| \le T$, $w \ge y \ge D$, $T \ge D \ge 2$.  In particular, for the value of $t$ in the definition of $\rho$
\begin{equation}
\Re \sum\limits_{D < p \le x} a_p^2 \chi(p) p^{-1-it} \le \log L_2 + O(1), \notag
\end{equation}
where $L_2$ denotes $\log(\log x/\log D)$.

For each positive $\lambda$
\begin{align}
\sum\limits_{D < p \le x} a_p^2\left( 1 - \Re \chi(p) p^{it} \right) p^{-1} & \le \lambda^2 \sum\limits_{\substack{D < p \le x \\ |a_p| \le \lambda}} g(p) \left( 1 - \Re \chi(p) p^{it} \right) p^{-1} + 2\lambda^{-2} \sum\limits_{\substack{D < p \le x \\ |a_p| > \lambda}} a_p^4 p^{-1} \notag \\
& \le \lambda^2 \rho + 4\lambda^{-2} \left( L_2 + O(1) \right) \notag
\end{align}
and
\begin{equation}
\rho \ge \lambda^{-2} \left( 1 - 4\lambda^{-2} \right)L_2 - \log L_2 + O\left(\lambda^{-4} + 1 \right). \notag
\end{equation}
The coefficient of $L_2$ is maximized by the choice $\lambda = 2\sqrt{2}$, whence $\rho \ge (1/16)L_2 - \log L_2 + O(1)$.

Application of Theorem \ref{JP_elliott_04_thm_02} with $c = 1/2$, $\beta = 1$ obtains for a typical exceptional sum the bound
\begin{equation}
\sum\limits_{n \le x} g(n) \chi_j(n) \ll \frac{x}{\log x} \prod\limits_{\substack{p \le D \\ (p,D)=1}} \left( 1 + \frac{g(p)}{p} \right) \left(\frac{\log D}{\log x} \right)^{1/48} \left(\log\left( \frac{\log x}{\log D} \right)\right)^{1/3}. \notag
\end{equation}

Bearing in mind the note preceding the statement of Theorem \ref{JP_elliott_04_thm_02}, the asserted result follows rapidly from applications of Lemmas \ref{JP_elliott_04_lem_02} and \ref{JP_elliott_04_lem_03}.

This completes Example \ref{JP_elliott_04_exa_04} and our consideration of uniformities.
\end{example}


\emph{Remark}.  Although without immediate uniformity in $D$, the argument of Example \ref{JP_elliott_04_exa_04} readily furnishes examples in the application of Theorem \ref{JP_elliott_04_thm_08}.

\emph{Closing comments}.  In this study $\mathbb{Q}^*$ denotes the multiplicative group of positive rationals.  Otherwise, notation is largely standard.  In particular, $\mathbb{C}^*$, $(\mathbb{Z}/D\mathbb{Z})^*$ denote the multiplicative group of the nonzero complex numbers and of the reduced residue classes $\modd{D}$, respectively.

Symbols for standard number theoretical functions such as the greatest integer function, $[ \, \cdot \, ]$, may be found in Hardy and Wright \cite{Hardy&Wright2008}.

Definition of the ad\`elic group $GL_2(A_F)$ and the derivation of automorphic forms from its representations may be found in Gelbart \cite{gelbart1975automorphicadele}.

A detailed study of Maass forms attached to representations of $GL_k(A_\mathbb{Q})$ may be found in Goldfeld \cite{Goldfeld2006}.  Shimura's analytic continuation of braided symmetric square $L$-functions, appraised for application to Maass forms, may be found in Chapter 7 of that volume, within a treatment of the Gelbart-Jacquet lift.

\bibliographystyle{amsplain}
\bibliography{MathBib}

First published by transmission September 4, 2013.  An occasional typographical oversight may remain.

\end{document}